\documentclass[twoside,10pt]{amsart}

\usepackage{enumerate}
\usepackage{graphics}
\usepackage{graphicx}
\usepackage{amsmath}
\usepackage{amssymb}
\usepackage{amstext}
\usepackage{times}
\usepackage[latin1]{inputenc}
\usepackage{latexsym}
\usepackage{amsthm}
\usepackage{amscd}

\def\R{\mathbb{R}}
\def\C{\mathbb{C}}

\def\N{\mathbb{N}}

\def\B{\mathbb{B}}
\def\D{\mathbb{D}}

\def\esf{\mathbb{S}}

\newcommand{\metri}[1]{ds_{#1}}

\newcommand{\cte}{{\rm const \:}}

\newcommand{\longui}{\operatorname{length}}
\newcommand{\dist}{\operatorname{dist}}

\newcommand{\wect}[3]{#3_{#1} , \ldots , #3_{#2}}

\newcommand{\Int}[1]{{\operatorname{Int}}\,{#1}}
\newcommand{\intc}{\operatorname{Int}}
\newcommand{\mJ}{\mathcal{J}}

\def\a{{\alpha}}
\def\t{{\theta}}
\def\T{{\Theta}}
\def\g{{\gamma}}

\def\l{{\lambda}}
\def\L{{\Lambda}}
\def\de{{\delta}}
\def\be{{\beta}}

\def\vp{{\varphi}}
\def\om{{\omega}}
\def\Om{{\Omega}}
\def\s{{\sigma}}
\def\ep{{\epsilon}}

\newtheorem{lemma}{Lemma}
\newtheorem{remark}{Remark}
\newtheorem{theorem}{Theorem}
\newtheorem{proposition}{Proposition}
\newtheorem{corollary}{Corollary}
\newtheorem{definition}{Definition}
\newtheorem{claim}{Claim}[section]
\newtheorem*{theoremintro}{Theorem}
\newtheorem*{theoremintroA}{Theorem A}
\newtheorem*{theoremintroB}{Theorem B}
\newtheorem*{theoremintroC}{Theorem C}

\numberwithin{equation}{section}

\begin{document}

\title{Density theorems for complete minimal surfaces in $\R^3$}
\author{Antonio Alarcón}
\author{Leonor Ferrer}
\author{Francisco Martín}

\thanks{This research is partially supported by MEC-FEDER Grant no. MTM2004 - 00160.}
\date{\today}

\address{Departamento de Geometría y Topología\hfill\break\indent  Universidad de Granada, \hfill\break\indent18071, Granada \hfill\break\indent Spain}

\email[A. Alarcón]{alarcon@ugr.es}
\email[L. Ferrer]{lferrer@ugr.es}
\email[F. Martín]{fmartin@ugr.es}

\begin{abstract} 
In this paper we have proved several approximation theorems for the family of minimal surfaces in $\R^3$ that imply, among other things, that complete minimal surfaces are dense in the space of all minimal surfaces endowed with the topology of $C^k$ convergence on compact sets, for any $k \in \N$. 

As a consequence of the above density result, we have been able to produce the first example of a complete proper minimal surface in $\R^3$ with uncountably many ends.

\noindent {\em 2000 Mathematics Subject Classification. Primary 53A10; Secondary 49Q05, 49Q10, 53C42.
}
\end{abstract}

\maketitle


\section{Introduction}

The conformal structure of a complete minimal surface in $\R^3$ influences many of its global properties. A complete (orientable) minimal surface  has an underlying complex structure that can be either {\em parabolic} or {\em hyperbolic} (the elliptic (compact) case is not possible for a minimal surface in Euclidean space.) Classically, a Riemann surface without boundary is called hyperbolic if it  carries  a nonconstant positive superharmonic function and parabolic if it is neither compact nor hyperbolic.

Until the 1980's, it was a general thought that complete minimal surfaces of hyperbolic type played a marginal role in the global theory of minimal surfaces. However, the techniques and methods developed to study the Calabi-Yau problem  have showed that these surfaces are present in some of the most interesting aspects of the theory. It is natural that the first examples of complete hyperbolic minimal surfaces appeared as counterexamples to the Calabi-Yau conjectures, which original statement was given in 1965 by E. Calabi \cite{calabi} (see also \cite{chern} and \cite{Yau}). This author conjectured that {\em ``a complete minimal hypersurface in $\R^n$ must be unbounded''}, even more, {\em ``a complete nonflat minimal hypersurface in $\R^n$ has an unbounded projection in every $(n-2)$-dimensional affine subspace''}. 

Both conjectures turned out to be false, at least in the immersed case. In 1980, L. P. Jorge and F. Xavier \cite{Jorg-Xavi} constructed complete nonflat minimal disks in an open slab of $\R^3$ giving a counterexample  to the second conjecture. An important progress came in 1996, when N. Nadirashvili \cite{Nadirashvili} constructed the first example of a complete bounded minimally immersed disk in $\R^3.$ Initially, Nadirashvili's work seemed to be the end point of a classical problem. However, the methods and ideas introduced by this author were the beginning of a significant development in the construction of complete hyperbolic minimal surfaces. So, it has been possible to find examples with more interesting topological and geometrical properties. At the same time, some non-existence theorems have imposed some limits to the theory. Three have been the main lines of study.

{\em Embeddedness} creates a dichotomy in the Calabi-Yau's question. T. Colding and W. P. Minicozzi \cite{Cold-Mini} have proved that a complete embedded minimal surface with finite topology in $\R^3$ must be properly embedded in $\R^3$. In particular it cannot be contained in a ball. Very recently, Colding-Minicozzi result has been generalized  in two different directions. On one hand W. H. Meeks III, J. Pérez and A. Ros \cite{mpr-2} have proved  that if $M$ is a complete embedded minimal surface in $\R^3$ with finite genus and a countable number of ends, then $M$ is properly embedded in $\R^3$. On the other hand, Meeks and Rosenberg  have obtained that if a complete embedded minimal surface $M$ has injectivity radius $I_M >0$,  then $M$ is  proper in space. This is a corollary of the minimal lamination closure theorem \cite{mr-lam}.
As a consequence of the above results, it has been conjectured by Meeks, Pérez and Ros that {\em ``if $M \subset \R^3$  is a complete embedded minimal surface with finite genus, then $M$ is proper''.}
We would like to mention that the conjecture seems to be false under the assumption of infinite genus, as Meeks is working in the existence of a complete embedded minimal surface with infinite genus which is contained in a half space \cite{meeks-com}.

The second line of work is related with the {\em properness} of the examples. Recall that an immersed submanifold of $\R^n$ is proper if the pre-image through the immersion of any compact subset of $\R^n$ is compact in the submanifold. It is clear from the definition that a proper minimal surface in $\R^3$ must be unbounded, so Nadirashvili's surfaces are not proper in $\R^3$. Much less obvious is that Nadirashvili's technique did not guarantee the immersion $f \colon \mathbb{D} \rightarrow \B$ was proper in the unit (open) ball $\B$ (here $\D$ stands for the unit disk in $\C$), where by proper we mean in this case that $f^{-1}(C)$ is compact for any $C\subset \B$ compact. Morales and the third author \cite{tran, MM-Convex, propi} introduced completely new ingredients in Nadirashvili's machinery and they proved that every convex domain (not necessarily bounded or smooth) admits a complete properly immersed minimal disk. These examples disproved a longstanding conjecture, which asserted that a complete minimal surface (without boundary) with finite topology and which is properly immersed in $\R^3$ should be parabolic. Recently \cite{MM-Convex2} they improved on their original techniques and were able to show that every bounded domain with $C^{2,\alpha}$-boundary admits a complete properly immersed minimal disk whose limit set is close to a prescribed simple closed curve on the  boundary of the domain. Similar methods of construction have been used by M. Tokuomaru in \cite{miki} to produce a complete minimal annulus properly immersed in the unit ball of $\R^3$. In contrast to these existence results for complete  properly immersed minimal disks in bounded domains, Meeks, Nadirashvili and the third author \cite{mhk} proved the existence of bounded open regions of $\R^3$ which do not admit complete properly immersed minimal surfaces with an annular end.  In particular, these domains do not contain a complete properly immersed minimal surface with finite topology.
\begin{figure}[htbp]
	\begin{center}
		\includegraphics[width=0.60\textwidth]{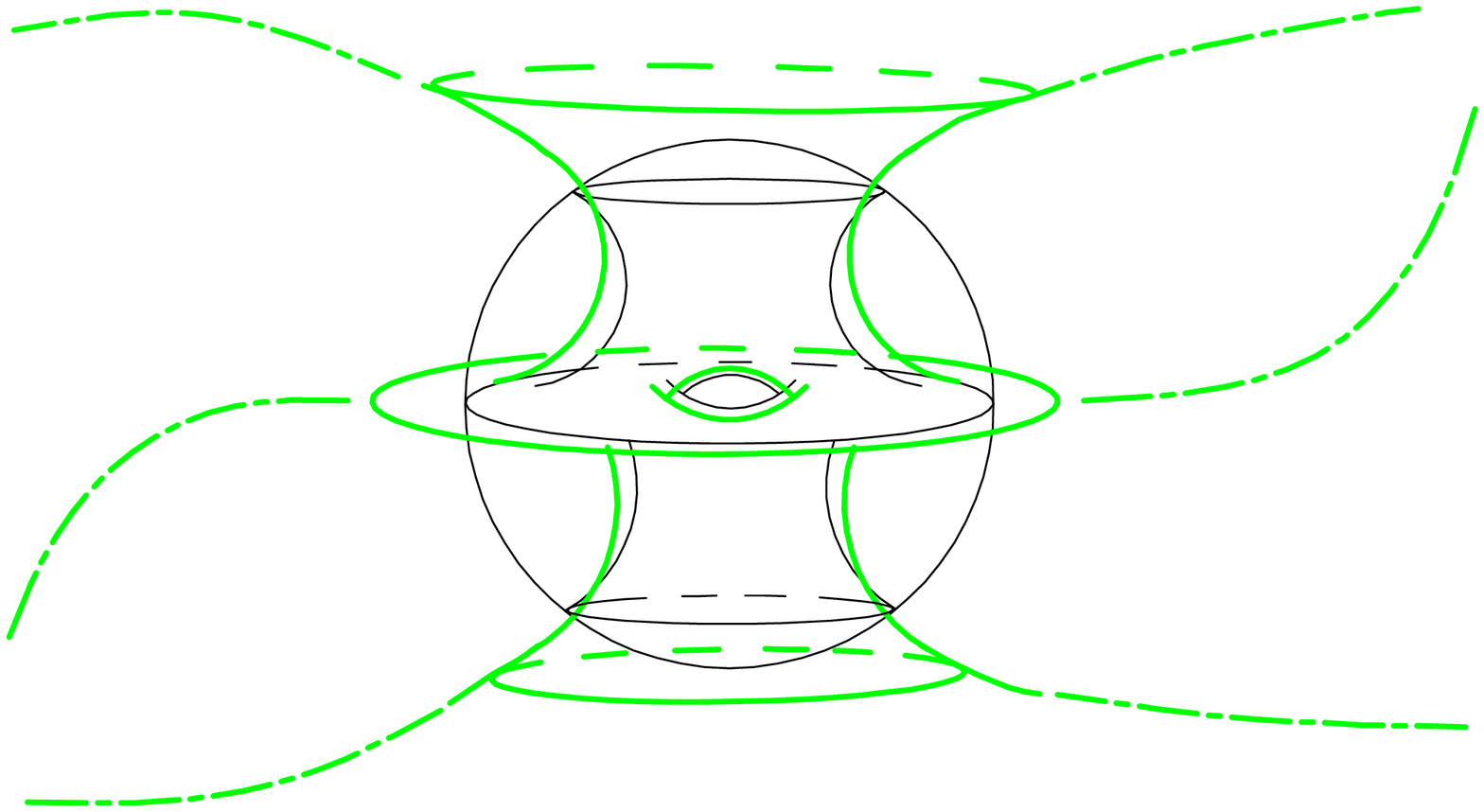}
	\end{center}
	\caption{} \label{fig:costa3}
\end{figure}

The other line of study for complete hyperbolic minimal surfaces in $\R^3$ has been the construction of {\em examples with nontrivial topology.} Nadirashvili's examples are simply connected. Thus, his mathematical machinery, which is based on Runge's theorem and López-Ros transformation, works without problems. López, Morales and the third author \cite{LMM-Handles, LMM-Nonorientable} introduce a third element in the construction: the Implicit Function Theorem, in order to produce Runge's functions that close also the periods when they are used as parameters in the López-Ros deformation.   

The aim of this paper is to join the second and third lines of work described in the above paragraphs in order to prove the following result (Section \ref{sec:teoremas}, Theorems \ref{teo-ultimo} and \ref{teo-B}).
\begin{theoremintroA}[\bf Density theorem] $\;$
Properly immersed, hyperbolic minimal surfaces of finite
topology are dense in the space of all properly immersed minimal surfaces in $\R^3$, endowed with the topology of smooth
convergence on compact sets.
\end{theoremintroA} 
Note that the best understood families of minimal surfaces in $\R^3$ (properly embedded, periodic, finite total curvature, finite type,...) are included in the statement of Theorem A. Furthermore, if we do not care about properness, then we can prove that:
\begin{quote} \em Complete (hyperbolic) minimal surfaces are dense in the space of minimal surfaces in $\R^3$ (without boundary) endowed with the topology of $C^k$ convergence on compact sets, for any $k \in \N$.\end{quote}

In the case of hyperbolic minimal surfaces we have an infinite number of linearly independent Jacobi fields. This is the key point in the proof of the above theorem. This enormous capability of deformation allows us to ``model'' a given compact piece of a hyperbolic minimal surface in order to approximate any other minimal surface with the same topological type (see Figure \ref{fig:costa3}). In particular, we can obtain the following existence result.
\begin{theoremintroB}
For any convex domain $D$ in $\R^3$ (not necessarily bounded or smooth) there exists a complete proper minimal immersion $\psi:M \rightarrow D$, where $M$ is an open  Riemann surface with arbitrary finite topology.
\end{theoremintroB}

One of the most interesting applications of our Density Theorem is the construction of the first example of a complete minimal surface properly immersed in $\R^3$ with an uncountable number of ends (Section \ref{sec:uncountable}). 
\begin{theoremintroC}
There exists a domain $\Omega \subset \C$ and a complete proper minimal immersion $\psi:\Omega \rightarrow \R^3$ which has uncountably many ends.
\end{theoremintroC}
The domain $\Omega$ is bounded in $\C$ and its set of ends contains a Cantor's set. We would like to emphasize that our technique can be also applied to construct complete proper minimal surfaces of genus $k$, $k\in \N$, and uncountably many ends. For the sake of simplicity, we have only exhibited in this paper the construction of a minimal surface of genus zero and uncountably many ends.

Once again, embeddedness establishes a dichotomy in the global theory of minimal surfaces. So, it is important to note that complete proper minimal surfaces in $\R^3$ with uncountably many ends cannot be embedded as a consequence of a result by Collin, Kusner, Meeks and Rosenberg \cite{CKMR}.
\vskip .3cm

\noindent {\em Acknowledgments.} We are indebted to W. H. Meeks III for valuable suggestions in the construction of minimal surfaces with an uncountable number of ends. We would also like to thank A. Ros for helpful criticisms of the paper.




\section{Preliminaries}
 \label{sec:pre}

This section is devoted to briefly summarize the notation and
results about Riemann surfaces, minimal surfaces,  and convex geometry that we will use in the paper.

\subsection{Riemann surfaces background}

Throughout the paper  $M'$ will denote a connected compact Riemann surface of genus $\s\in\N \cup \{0\}.$ 

Let $M$ be a domain in $M'$ and assume that $M$ carries a Riemannian metric $ds^2$. Given a subset $W \subset M$, we define:
\begin{itemize}
\item $\dist_{(W,ds)}(p,q)=\inf \{\longui(\alpha, ds) \: | \: \alpha:[0,1]\rightarrow W, \; \alpha(0)=p,\alpha(1)=q \}$, for any $p,q\in W$;
\item $\dist_{(W,ds)}(T_1,T_2)=\inf \{\dist_{(W,ds)}(p,q) \;|\;p \in T_1, \;q \in T_2 \}$, for any $T_1, T_2 \subset  W$;
\item $\text{diam}_{ds}(W)=\sup \{\dist_{(W,ds)}(p,q) \;|\;p,q \in W \}$. 
\end{itemize}

For $\textsc{e}\in\N$, consider $\D_1,\ldots,\D_\textsc{e}\subset M'$ open disks so that $\{\g_i:=\partial \D_i\}_{i=1}^{\textsc{e}}$ are piecewise smooth Jordan curves and $\overline{\D}_i\cap \overline{\D}_j=\emptyset$ for all $i\neq j$.

\begin{definition}\label{multicycle}
Each curve $\g_i$ will be called a cycle on $M'$ and the family $\mathcal{ J}=\{\g_1,\ldots,\g_\textsc{e}\}$ will be called a multicycle on $M'$. We denote by $\intc(\g_i)$ the disk $\D_i$, for $i=1,\ldots, \textsc{e}.$ We also define $M(\mathcal{ J})=M'-\cup_{i=1}^\textsc{e} \overline{\intc(\g_i)}$. Notice that $M(\mathcal{ J})$ is always connected.
\end{definition}
Given $\mathcal{ J}=\{\g_1,\ldots,\g_\textsc{e}\}$ and $\mathcal{ J}'=\{\g_1',\ldots,\g_\textsc{e}'\}$ two multicycles in $M'$ we write $\mathcal{ J}'< \mathcal{ J}$ if $\overline{\intc(\g_i)} \subset \intc (\g_i')$ for $i=1,\ldots, \textsc{e}.$ Observe that this implies $\overline{M(\mathcal{ J}')} \subset M(\mathcal{ J})$.

Let $\mathcal{ J}=\{\g_1,\ldots,\g_\textsc{e}\}$ be a multicycle and assume that $\overline{M(\mJ)}\subset M$, where the Riemannian metric $ds^2$ is defined. If $\ep>0$ is small enough, we can consider the multicycle $\mathcal{ J}^\ep=\{\g_1^\ep,\ldots,\g_\textsc{e}^\ep\}$, where by $\g_i^\ep$ we mean the cycle satisfying $\overline{\intc(\g_i)} \subset \intc (\g_i^\ep)$ and $\dist_{(M,ds)}(q,\g_i)=\ep$ for all $q\in \g_i^\ep$ and $i=1,\ldots, \textsc{e}.$ Similarly, we can define $\mathcal{ J}^{-\ep}=\{\g_1^{-\ep},\ldots,\g_\textsc{e}^{-\ep}\}$, where $\g_i^{-\ep}$ now means the cycle satisfying $\overline{\intc(\g_i^{-\ep})} \subset \intc (\g_i)$ and $\dist_{(M,ds)}(q,\g_i)=\ep$ for all $q\in \g_i^{-\ep}$ and $i=1,\ldots, \textsc{e}.$

Given a Riemann surface with boundary $N \subset M'$, we will say that a function, or a 1-form, is harmonic, holomorphic, meromorphic, ... on $\overline N$, if it is harmonic, holomorphic, meromorphic, ... on a domain containing $\overline N$.

\subsection{Minimal surfaces background}
The theory of complete minimal surfaces is closely related to the theory of Riemann surfaces. This is due to the fact that any such surface is given by a triple $\Phi=(\Phi_1, \Phi_2, \Phi_3)$ of holomorphic 1-forms defined on some Riemann surface such that
\begin{equation} \label{eq:conforme}
\Phi_1^2+\Phi_2^2+\Phi_3^2=0;
\end{equation}
\begin{equation} \label{eq:bilbao}
\|\Phi_1\|^2+\|\Phi_2\|^2+\|\Phi_3\|^2 \neq 0;
\end{equation}
and all periods of the $\Phi_j$ are purely imaginary, here we consider $\Phi_i$ to be a holomorphic function times $dz$ in a local parameter $z$. Then the minimal immersion $X:M \rightarrow \R^3$ can be parameterized by $z \mapsto \mbox{Re} \int^z \Phi.$
The above triple is called the Weierstrass representation of the immersion $X$. Usually, the first requirement (\ref{eq:conforme}) (which ensures the conformality of $X$) is guaranteed by introducing the formulas:
$$\Phi_1 =\frac12 \left( 1-g^2\right) \, \eta, \quad \Phi_2 =\frac{\rm i}2 \left( 1+g^2\right) \, \eta, \quad \Phi_3= g \, \eta, $$
with a meromorphic function $g$ (the stereographic projection of the Gauss map) and a holomorphic 1-form $\eta$. The metric of $X$ can be expressed as
\begin{equation} \label{eq:metric}
{\metri X}^2=\tfrac12 \|\Phi\|^2=\left(\tfrac12\left(1+|g|^2\right) \|\eta \|\right)^2.
\end{equation}
Throughout the paper, we will use several orthonormal bases of $\R^3$. Given $X:\Omega \rightarrow \R^3$ a minimal immersion and $S$ an orthonormal basis, we will write the Weierstrass data of $X$ in the basis $S$ as $${\Phi}_{(X,S)}=(\Phi_{(1,S)},\Phi_{(2,S)},\Phi_{(3,S)}),  \quad g_{(X,S)}, \quad \eta_{(X,S)}.$$
Similarly, given $v \in \R^3$, we will let $v_{(k,S) }$ denote the $k$-th coordinate of $v$ in $S$. The first two coordinates of $v$ in this basis will be represented by $v_{(*,S)}=\left(v_{(1,S)},v_{(2, S)} \right)$.

Given a curve $\alpha$ in $M$, by $\longui(\alpha, X)$ we
mean $\longui(\alpha, \metri X)$. Similarly, given a subset $W \subset M$, we write:
\begin{itemize}
\item $\dist_{(W,X)}(p,q)=\dist_{(W,\metri X)}(p,q)$, for any $p,q\in W$;
\item $\dist_{(W,X)}(T_1,T_2)=\dist_{(W,\metri X)}(T_1,T_2)$, for any $T_1, T_2 \subset  W$;
\item $\text{diam}_{X}(W)=\text{diam}_{\metri X}(W)$. 
\end{itemize}

\subsubsection{The López-Ros transformation}
The proof of Lemmas \ref{AdaptaRunge-2}, \ref{properness} and \ref{nadira} exploits what has come to be called the López-Ros transformation. If $M$ is a Riemann surface and $(g,\eta)$ are the Weierstrass data of a minimal immersion $X:M \rightarrow \R^3$, we define on $M$ the data
\begin{equation} \label{eq:mambru}
\widetilde g= \frac{g}{h}, \qquad \widetilde \eta= \eta \cdot h,
\end{equation}
where $h:M \rightarrow \C$ is a holomorphic function without
zeros. If the periods of this new Weierstrass representation are purely imaginary, then it defines a minimal immersion $\widetilde X: M \rightarrow \R^3$. This method provides us with a powerful and natural tool for deforming minimal surfaces. From our point of view, the most important property of the resulting surface is that the {\em third coordinate function is preserved.}  Note that the intrinsic metric is given by (\ref{eq:metric}) as
 \begin{equation} \label{eq:metric2}
{\metri{\widetilde X}}^2=\left(\tfrac12\left(|h|+\frac{|g|^2}{|h|}\right)\, \|\eta \| \right)^2.
\end{equation}
This means that we can increase the intrinsic distance in a prescribed compact of $ M$, by using  suitable functions $h$. These functions will be provided by Lemma \ref{AdaptaRunge-2} that can be consider a Runge's type theorem.
\subsection{Background on convex geometry} 
Convex geometry is a classical subject with a large literature. To make this article self-contained, we will describe the concepts and results we will need. 
A convex, compact set of $\R^n$ with nonempty interior is called {\em a convex body}. A theorem of H. Minkowski (cf. \cite{mingorebulgo}) states that every convex body $C$ in $\R^n$ can be approximated (in terms of Hausdorff metric) by a sequence $C_k$ of `analytic' convex bodies. Recall that the Hausdorff distance between two nonempty compact subsets of $\R^n$, $C$ and $D$, is given by:
$$\delta^H(C,D)=\max \left\{ \sup_{x \in C} \inf_{y \in D} \|x-y \|, \sup_{y \in D} \inf_{x \in C} \|x-y \| \right\}.$$
\begin{theoremintro}[Minkowski] \label{th:minko}
Let $C$ be a convex body in $\R^n$. Then there exists a sequence $\{C_k \}$ of convex bodies with the following properties
\begin{enumerate}[\rm 1.]
\item $C_k \searrow C$ in terms of the Hausdorff metric;
\item $\partial C_k$ is an analytic $(n-1)$-dimensional manifold;
\item The curvatures of $\partial C_k$ never vanish.
\end{enumerate}
\end{theoremintro}
A modern proof of this result can be found in \cite[\S 3]{meeksyau}. 

Given $E$ a bounded regular convex domain of $\R^3$ and $p \in
\partial E$, we will let $\kappa_2(p) \geq \kappa_1(p) \geq 0$ denote the 
principal curvatures of $\partial E$ at $p$ (associated to the inward
pointing unit normal.) Moreover, we write: $$\kappa_1(\partial E) := \mbox{min} \{ \kappa_1(p) \: | \: p \in \partial E \}. $$ 
If we consider $\mathcal{N}: \partial E \rightarrow \esf^2$ the outward
pointing unit normal or Gauss map of $\partial E$, then there exists a constant $a>0$ (depending on $E$) such that $\partial E_t=\{p+ t\cdot \mathcal{N}(p) \; | \; p \in \partial E\}$ is a regular (convex) surface $\forall t \in [-a, +\infty[$. We label $E_t$ as the convex domain bounded by $\partial E_t$. The normal projection to $\partial E$ is represented as 
$$\begin{array}{rccl} \mathcal{P}_E: &\R^3- E_{-a}& \longrightarrow & \partial E \\ \quad & p+t \cdot \mathcal{N}(p) & \longmapsto & p \; .\end{array}$$
Finally, we define the `extended' Gauss map $\mathcal{N}_E : \R^3 - E_{-a} \longrightarrow \esf^2$ as $\mathcal{N}_E(x)=\mathcal{N}(\mathcal{P}_E(x)).$




\section{A Runge's type Lemma} \label{sec:runge}
As we mentioned in the introduction, this section contains a Runge type theorem on Riemann surfaces. It will be crucial in the prove of the main theorems. 

\begin{lemma}\label{AdaptaRunge-2}
Let $\mathcal{ J}$ be a multicycle of $M'$ and $F:\overline{M(\mathcal{ J})}\to\R^3$ a conformal minimal immersion with Weierstrass data $(g,\Phi_3)$. Consider $K_1$ and $K_2$ two disjoint compact sets in $M(\mathcal{ J})$ and $\Delta\subset M'$ satisfying:
\begin{itemize}
\item There exists a basis of the homology of $M(\mathcal{J})$ contained in $K_2;$
\item $\overline{\Delta}\subset M'-(K_1\cup K_2);$ 
\item $\Delta$ has a point in each connected component of $M'-(K_1\cup K_2).$
\end{itemize}
Then, for any $m\in\N$ and for any $t>0$ there exists $H:\overline{M(\mathcal{ J})}-\Delta\to\C$ a holomorphic function without zeros, such that
\begin{enumerate}[\rm (L\ref{AdaptaRunge-2}.a)]
\item $\left|H-t\right|<1/m$ in $K_1$;
\item $\left|H-1\right|<1/m$ in $K_2$;
\item The minimal immersion $\widetilde{F}:\overline{M(\mathcal{J})}-\Delta\to\R^3$ with Weierstrass data $(g/H,\Phi_3)$ is well-defined.
\end{enumerate}
\end{lemma}


In order to prove Lemma \ref{AdaptaRunge-2}, we have to introduce some terminology and prove several claims. We define $\varrho=2\s+\textsc{e}-1$ (recall that $\s$ is the genus of the compact surface $M'.$) Thus, let $\mathcal{ B}=\{\aleph_1,\ldots,\aleph_\varrho\}$ be a basis of the homology of $M(\mJ)$ contained in $K_2,$ and denote by $\mathcal{ H}$ the complex vector space of the holomorphic 1-forms on $\overline{M(\mJ)}.$

\begin{claim}\label{1Runge}
Consider $(\wect{1}{\varrho}{a})\in\C^{\varrho}-\{(0,\dots,0)\}$ and $c=\sum_{j=1}^{\varrho}a_j\aleph_j$. Then there exists $\tau\in \mathcal{ H}$ with $\int_c\tau = 1$.
\end{claim}

\begin{proof}
The first holomorphic De Rham cohomology group, $H_{\text{hol}}^1(\overline{M(\mJ)})$ is a complex vector space of dimension $\varrho$ (see \cite[Chapter III.5]{Farkas}). Thus, the map $I:H_{\text{hol}}^1(\overline{M(\mJ)})\to\C^{\varrho}$ given by $I([\psi])=\left(\int_{\aleph_1} \psi,\ldots,\int_{\aleph_\varrho}\psi\right),$ is a linear isomorphism. Observe $I$ is well-defined from the fact that the type of an exact 1-form in $H_{\text{hol}}^1(\overline{M(\mJ)})$ is zero. Therefore, there exists $[\psi]\in H_{\text{hol}}^1(\overline{M(\mJ)})$ such that $I([\psi])\notin \{(\wect{1}{\varrho}{z})\;|\; \sum_{j=1}^\varrho a_jz_j=0\}$. Therefore, we can choose $\tau\in [\psi]$ with $\int_c\tau =1.$
\end{proof}


\begin{claim}\label{2Runge}
Consider $\tau\in \mathcal{ H}$ and $P\in \overline{M(\mJ)}$. Then, there exists a holomorphic function $A:\overline{M(\mJ)}\to\C$ such that $\left(\left.(\tau+dA)\right|_{\overline{M(\mathcal{J})}}\right)_0\geq \left(\left.\tau\ \right|_{\overline{M(\mathcal{ J})}}\right)_0\cdot P$, where $(\cdot)_0$ denotes the divisor of zeros.
\end{claim}

\begin{proof}
Suppose $(\tau)_0=Q_1\cdots Q_k P^n,$ with $P\neq Q_i$ $\forall\, i=1,\ldots,k,$ and assume that there exists a holomorphic function $\upsilon:\overline{M(\mJ)}\to\C$ satisfying
\begin{enumerate}[1)]
\item $\upsilon(P)\neq \upsilon(Q_i),$ $\forall\,i=1,\ldots,k;$
\item $P$ is not a ramification point of $\upsilon.$
\end{enumerate}
Consider the function $J:\overline{M(\mJ)}\to\C$ given by $J=(\upsilon-\upsilon(P))^{n+1}\prod_{i=1}^k(\upsilon-\upsilon(Q_i))^2.$ Therefore, $\left(\left. dJ \right|_{\overline{M(\mathcal{ J})}}\right)_0\geq \left(\left.\tau\right|_{\overline{M(\mathcal{ J})}}\right)_0$ and the order of $P$ as zero of $dJ$ and $\tau$ is the same, so, there exists $\l\in\C$ such that $A=\l J$ solves the claim.

Now, we are checking that there exists a such function $\upsilon$ satisfying items 1) and 2). A Runge's type theorem (see \cite[Theorem 10]{Royden}) guarantees the existence of a holomorphic function $\upsilon_1:\overline{M(\mJ)}\to\C$ fulfilling item 1). On the other hand, given $(U,z)$ a conformal coordinate chart around $P$ and $m\in\N,$ the same theorem provides us of a holomorphic function $h_m:\overline{M(\mJ)}\to\C$ with $|h_m(z)-z|<1/m$ for $z  \in U$. Hence, $\{h_m\}_{m\in\N}\to z$ and therefore $\{dh_m\}_{m\in\N}\to dz.$ Taking into account that $P$ is not a ramification point of $z$, we conclude that there exists $m\in\N$ large enough so that $P$ is not a ramification point of $\upsilon_2:=h_m.$ Finally we choose $\upsilon$ as a appropriate linear combination of $\upsilon_1$ and $\upsilon_2.$
\end{proof}


\begin{claim}\label{3Runge}
Let $\mathcal{O}(\overline{M(\mathcal{ J})})$ be the real vector space of the holomorphic functions on $\overline{M(\mathcal{ J})}$. Then the linear map $\mathcal{F}:\mathcal{O}(\overline{M(\mathcal{ J})})\to\R^{2\varrho}$ given by
$$\mathcal{F}(\vp)=\left( \text{\rm Re}\left[ \int_{\aleph_j} \vp\,\Phi_3\left( \frac{1}{g}+g \right) \right]_{j=1,\dots,\varrho}, \text{\rm Im}\left[ \int_{\aleph_j} \vp\,\Phi_3\left( \frac{1}{g}-g \right) \right]_{j=1,\dots,\varrho} \right)$$
is onto.
\end{claim}

\begin{proof}
Suppose $\mathcal{F}$ is not onto. Therefore, there exists $(\wect{1}{2\varrho}{\mu})\in\R^{2\varrho}-\{(0,\dots,0)\}$ such that $\mathcal{F} \left(\mathcal{O}(\overline{M(\mathcal{ J})})\right)$ $\subset$ $\{(\wect{1}{2\varrho}{x})\;|\; \sum_{j=1}^{2\varrho} \mu_jx_j=0\}.$ In other words:
\begin{equation}\label{R1}
\text{\rm Re}\left[\sum_{j=1}^{\varrho}\left( u_j\int_{\aleph_j}\frac{\vp}{g}\Phi_3+\overline{u_j}\int_{\aleph_j} \vp g\Phi_3\right) \right] =0 \; ,\quad \forall\, \vp\in \mathcal{O}(\overline{M(\mathcal{ J})})\;,
\end{equation}
where $u_j=\mu_j-i\mu_{j+\varrho},$ $j=1,\dots,\varrho$.

Now, Claims \ref{1Runge} and \ref{2Runge} guarantee the existence of a differential $\tau\in\mathcal{ H}$ satisfying
\begin{itemize}
\item $(\tau)_0\geq \left( \left. \left(\frac{1}{g}\Phi_3\right)\right|_{\overline{M(\mathcal{ J})}}\right)_0^2\cdot \left( \left. \left(g\, dg\right)\right|_{\overline{M(\mathcal{ J})}}\right)_0;$
\item ${\rm Re}\left[\sum_{j=1}^{\varrho} \overline{u_j}\int_{\aleph_j}\tau\right]=1.$
\end{itemize}
Therefore, if we define $w:=\frac{\tau}{2 g\, dg}$, then $\vp:=\frac{g\, d\,w}{\Phi_3}\in \mathcal{O}(\overline{M(\mathcal{ J})})$. Hence, integrating (\ref{R1}) by parts, we obtain
$$
\text{\rm Re}\left[\sum_{j=1}^{\varrho}\overline{u_j} \int_{\aleph_j} \vp g\Phi_3\right] =-\text{\rm Re}\left[\sum_{j=1}^{\varrho} \overline{u_j}\int_{\aleph_j}\tau\right]=0\;,
$$ which is absurd. This proves the claim.
\end{proof}


Using the above claim we obtain the existence of $\{\wect{1}{2\varrho}{\vp}\}\subset \mathcal{O}(\overline{M(\mathcal{ J})})$ such that $\{\mathcal{F}(\vp_1),$ $\dots,$ $\mathcal{F}(\vp_{2\varrho})\}$ \label{orejon} are linearly independent. Fixed $m_0\in\N$, without loss of generality, we can assume
\begin{equation}\label{R2}
\left| \sum_{i=1}^{2\varrho}x_i \vp_i(p) \right|<\frac{1}{m_0}\;,
\end{equation}
$\forall\, x=(\wect{1}{2\varrho}{x})\in\R^{2\varrho}$ with $\left\|x\right\|_\infty<1$, $\forall\, p\in\overline{M(\mathcal{ J})}$.


\subsection{Proof of Lemma \ref{AdaptaRunge-2}}

Given $n\in \N,$ we apply a Runge-type theorem on $M',$ see \cite[Theorem 10]{Royden}, and obtain a holomorphic function $\vartheta_n:\overline{M(\mathcal{ J})}-\Delta\to\C$ such that
\begin{equation*}\label{phi-cte}
\begin{cases}
|\vartheta_n-n\,\log(t)|<1/n & \text{in } K_1 \; ,\\
|\vartheta_n|<1/n & \text{in } K_2\;.
\end{cases}
\end{equation*}
Now, for $\T=(\l_0,\ldots,\l_{2\varrho})\in\R^{2\varrho+1},$ we consider the map $h^{\T,n}:\overline{M(\mathcal{ J})}-\Delta\to\C$ given by
$$
h^{\T,n}(p)=\exp\left[ \l_0\, \vartheta_n(p)+\sum_{j=1}^{2\varrho}\l_j\vp_j(p) \right]\;.
$$
Label $g^{\Theta,n}=g/h^{\T,n}$ and $\Phi_3^{\Theta,n}=\Phi_3$. Clearly, we have that $\vartheta_n$ converges uniformly on $K_2$ to $\vartheta_\infty\equiv 0$. So, for $\T=(\wect{0}{2\varrho}{\l})\in \R^{2\varrho+1}$ we also define on $K_2$ the Weierstrass data $g^{\T,\infty}=g/h^{\T,\infty}$ and $\Phi_3^{\T,\infty}=\Phi_3$, where $h^{\T,\infty}:K_2\to\C$ is given by
\begin{equation}\label{h-infi}
h^{\T,\infty}(p)=\text{exp}\left[\sum_{j=1}^{2\varrho}\l_j\, \vp_j(p)\right]\;.
\end{equation}

Note that the third coordinate of all these Weierstrass representations has no real periods, but the period problems of the two first ones coordinates are not solved. In order to solve these problems we define, $\forall\, n\in\N\cup\{\infty\}$, the map $\mathcal{ P}_n:\R^{2\varrho+1}\to\R^{2\varrho}$ given by 
\begin{equation}\label{Pn}
\mathcal{ P}_n(\T)=\left( {\rm Re}\left[\int_{\aleph_j}\Phi_1^{\T,n}\right]_{j=1,\dots,\varrho}\, ,\, {\rm Re}\left[ \int_{\aleph_j}\Phi_2^{\T,n}\right]_{j=1,\dots,\varrho}\right)\;.
\end{equation}

Since $\mathcal{F}$ is a well-defined immersion, then we have $\mathcal{ P}_n(0,\dots,0)=0$, $\forall\, n\in\N\cup\{\infty\}$. Moreover, it is not hard to check that 
$$
[\text{Jac}_{\wect{1}{2\varrho}{\l}}(\mathcal{ P}_n)](0,\dots,0)=\text{det}(\mathcal{F}(\vp_1),\dots,\mathcal{F}(\vp_{2\varrho}))\neq 0,\quad \forall\, n\in\N\cup\{\infty\}\;.
$$
Labeling $\overline{B}(0,r)=\{\L\in\R^{2\varrho}\;|\; \|\L\|\leq r\}$, we can find $\xi>0$ and $0<r<1$ such that the Jacobian operator $[\text{Jac}_{\wect{1}{2\varrho}{\l}}(\mathcal{ P}_\infty)]\big|_{[-\xi,\xi]\times\overline{B}(0,r)}\neq 0$ and $\mathcal{ P}_\infty(0,\cdot)\big|_{\overline{B}(0,r)}$ is injective.

As $\{\vartheta_n\}_{n\in\N}$ uniformly converges to $\vartheta_\infty\equiv 0$ on $K_2$ and $\aleph_i\subset K_2,$ $\forall\, i=1,\dots,\varrho$, then it is not hard to see that $\{\text{Jac}_{\wect{1}{2\varrho}{\l}}(\mathcal{ P}_n)\}_{n\in\N}$ uniformly converges to $\text{Jac}_{\wect{1}{2\varrho}{\l}}(\mathcal{ P}_\infty)$ on $[0,\xi]\times\overline{B}(0,r)$. Therefore, there exists $n_0\in\N$ satisfying that $\forall\, n\geq n_0,$ $\exists\ \xi_n>0$ such that $[\text{Jac}_{\wect{1}{2\varrho}{\l}}(\mathcal{ P}_n)](\l_0,\L)\neq 0$, $\forall\, (\l_0,\L)\in [-\xi_n,\xi]\times\overline{B}(0,r).$
Now, we are able to apply the Implicit Function Theorem to the map $\mathcal{ P}_n$ at $(0,\dots,0)\in [-\xi_n,\xi]\times\overline{B}(0,r)$ and obtain a smooth function $L_n:I_n\to\R^{2\varrho}$, satisfying $\mathcal{ P}_n(\l_0,L_n(\l_0))=0$, $\forall\, \l_0\in I_n$, where $I_n$ is an open interval containing 0 and maximal, in the sense that $L_n$ can not be regularly extended beyond $I_n$.

\begin{claim}\label{4Runge}
There exist $\ep_0>0$ and $n_0\in\N$ such that the function $L_n:[0,\ep_0]\to\overline{B}(0,r)$ is well-defined for all $n\geq n_0$.
\end{claim}

The proof of Claim \ref{4Runge} is a standard argument of classical analysis that can be found in \cite{LMM-Handles}.

Take $n\geq n_0$ large enough so that $1/n\leq\ep_0$  and label $(\wect{1}{2\varrho}{\l^{n}})=L_n(1/n)$. If $m_0$ in (\ref{R2}) and  $n\geq n_0$ are sufficiently large, the function
$$
H(p)=\text{exp}\left[\frac{1}{n}\vartheta_n(p)+\sum_{j=1}^{2\varrho}\l_j^{n}\, \vp_j(p)\right]
$$
satisfies (L1.a) and (L1.b). As the period function $\mathcal{ P}_n$ vanishes at $\T_n=(1/n,\wect{1}{2\varrho}{\l^{n}})$, then the minimal immersion $\widetilde{F}$ with Weierstrass data given by $(g/H,\Phi_3)$ is well-defined. Hence, the function $H$ also satisfies (L1.c). This completes the proof of Lemma \ref{AdaptaRunge-2}.





\section{Properness Lemma}
This lemma asserts that a compact minimal surface  whose boundary is close to the boundary of a convex $E$ can be `elongated' in such a way that the  boundary of the new surface achieves the boundary of a bigger convex $E'$. However, the above procedure does not change the topological type of the minimal surface. If $E$ is strictly convex we are able to obtain some extra information about the resulting surface that will be necessary  in proving Theorem \ref{teo-B} (see Remark \ref{k1>0}.) 

\begin{lemma}\label{properness}
Let $E$ and $E'$ be two  bounded regular convex domains in $\R^3$, with $0 \in E\subset \overline{E}\subset E'$. Consider $\mathcal{ J}'< \mathcal{ J}_0$ multicycles in $M'$ and $X:\overline{M(\mathcal{ J}_0)}\to\R^3$ a conformal minimal immersion  satisfying $X(p_0)=0$ for a given point $p_0 \in M(\mathcal{J}')$, and
\begin{equation}\label{(3)}
X( \, \overline{M(\mathcal{ J}_0)}-M(\mathcal{ J}') \,)\subset E'- \overline{E}\;.
\end{equation}
Finally, consider $b_2>0$ such that $E'_{-b_2}$ and $E_{-2  b_2}$ exist. Then, for any $b_1>0$ there exist a multicycle $\mathcal{ J}$ and a conformal minimal immersion $Y:\overline{M(\mathcal{ J})}\to\R^3$ satisfying $Y(p_0)=0$ and:
\begin{enumerate}[\rm (L\ref{properness}.a)]
\item $\mathcal{ J}'<\mathcal{ J}<\mathcal{ J}_0$;
\item $\left\|Y(p)-X(p)\right\|< b_1$, $\forall\, p\in\overline{M(\mathcal{ J}')}$;
\item $Y(\mathcal{ J})\subset E'-E'_{-b_2}$;
\item $Y(\overline{M(\mathcal{ J})}- M(\mathcal{ J}'))\subset\R^3- E_{-2 b_2}.$
\end{enumerate}
\end{lemma}


\subsection{Proof of Lemma \ref{properness}}

Let $\omega$ be a meromorphic differential on $M'$ so that $\omega$ has neither zeroes nor poles on $\overline{M(\mathcal{J}_0)}$. Then, it is well known that $ds^2:=\| \omega\|^2$ is a flat Riemannian metric on $\overline{M(\mathcal{J}_0)}$. 
\begin{remark}[\bf developing map] \label{re:developing}
Fixed a point $q \in \overline{M(\mathcal{J}_0)}$ the multivalued map given by:
\begin{equation} \label{eq:developing}
\mathsf{f}(p):= \int_q^p \omega,
\end{equation}
is called the developing map of $\|\omega\|^2$. It is known that $\|\omega \|^2=\mathsf{f}^* ds_0^2$, where $ds_0^2$ represents the Euclidean metric of $\C$. In particular, $\mathsf{f}$ can be seen as a local isometry. 
\end{remark}
Given $n\in\N$ we define an order relation in the set $I\equiv \{1,\ldots,n\}\times\{1,\ldots,\textsc{e}\}.$ We say $(j,l)>(i,k)$ if one of the two following situations occurs: $l=k$ and $j>i$ or $l>k.$ Moreover given $p\in M(\mathcal{J}_0)$ and $r>0$, we denote $D(p,r)=\{q\in M(\mathcal{J}_0)\;|\;\dist_{(M(\mathcal{J}_0),ds)}(p,q)<r\}.$ We also define two important constants that are chosen as follows:
\begin{itemize}
\item $\mu=\max\{\dist_{\R^3}(x,\partial E)\;\;|\;\;x\in E'\};$
\item $\ep_0>0$  which will only depend on the data of Lemma \ref{properness} (i.e., $X$, $\mathcal{J}_0$, $\mathcal{J}'$, $E$, $E'$, $b_1$, and $b_2$.) This positive constant will be determined later and it must be small enough to satisfy several inequalities appearing in this section.
\end{itemize}


\subsubsection{The first deformation}

\begin{claim}\label{puntos-p}
There exist a multicycle $\mathcal{ J}_1$  such that $\mathcal{ J}'<\mathcal{ J}_1<\mathcal{ J}_0$, and a set of points $\{ p_i^k\;|\;(i,k)\in I\}$ included in $M(\mathcal{ J}_1)- \overline{M(\mathcal{J}')}$, satisfying the following properties:
\begin{enumerate}[\rm 1)]
\item For any $k$, there exists a cycle $\g_k$ passing trough $\{p_1^k,\ldots,p_n^k\}$ (orderly) and contained in $M(\mathcal{ J}_1)-\overline{M(\mathcal{J}')};$ 

\item $\mathcal{ J}_2=\{\g_1,\ldots,\g_\textsc{e}\}$ is a multicycle with $\overline{M(\mathcal{J}')}\subset M(\mathcal{ J}_2);$

\item There exist open disks $B^{i,k}\subset M(\mathcal{ J}_1)-\overline{M(\mathcal{J}')}$ satisfying  $p_i^k,p_{i+1}^k\in B^{i,k}$, and such that (we adopt the convention $p_{n+1}^k=p_1^k$) 
\begin{equation}\label{Bik-peque}
\|X(p)-X(p')\|<\epsilon_0\;, \quad \forall\, p, p' \in B^{i,k},\quad \forall\,(i,k)\in I\;;
\end{equation}

\item For any $(i,k)\in I,$ there exists an orthonormal basis of $\R^3$, $S_i^{k}=\{e_1^{i,k},e_2^{i,k},e_3^{i,k}\}$, with $e_1^{i,k}=\mathcal{N}_E(X(p_i^k))$, and satisfying
\begin{equation}\label{basesS-1}
\left\|e_j^{i,k}-e_j^{i+1,k}\right\|<\frac{\ep_0}{3\mu}\;,\quad \forall j \in \{1,2,3\}\quad (e_j^{n+1,k}:=e_j^{1,k})\;,
\end{equation}
and
\begin{equation}\label{basesS-2}
f_{(X,S_i^{k})}(p_i^k)\neq 0\;, \quad \text{where} \quad f_{(X,S_i^{k})}:= \frac{\eta_{(X,S_i^{k})}}{\omega}\;;
\end{equation}

\item For  each $(i,k)\in I,$ there exist a complex constant $\theta_i^k$  which satisfies $|\theta_i^k|=1$, ${\rm Im}\, \theta_i^k\not=0$, and
\begin{equation} \label{zetas}
\left|\theta_i^k\frac{ f_{(X,S_i^k)}(p_i^k)}{|f_{(X,S_i^k)}(p_i^k)|}- 1 \right| < \frac{\ep_0}{3\mu}\;.
\end{equation}
\end{enumerate}
\end{claim}
\begin{proof}
Since $\mathcal{J}'$ is a set of piecewise regular curves, then we know that $\mathcal{N}_E(X(\mathcal{J}'))$ omits an open set $U$ of $\esf^2$. Hence, we can get a multicycle $\mathcal{ J}_1$ with $\mathcal{ J}'<\mathcal{ J}_1<\mathcal{ J}_0$ and $\mathcal{N}_E(X(M(\mathcal{ J}_1) - M(\mathcal{ J}'))) \subset \esf^2 - U$. Let $V_1$ and $V_2$ be a smooth orthonormal basis of tangent vector fields on $\esf^2 - U$. Then, we define
$ \xi_1(p)=\mathcal{N}_E(X(p)),$ $ \xi_2(p)=V_1\left(\mathcal{N}_E(X(p))\right)$ and $\xi_3(p)=V_2\left(\mathcal{N}_E(X(p))\right)$, $\forall p\in M(\mathcal{ J}_1)-\overline{M(\mathcal{ J}')}$.

If $n$ is large enough, because of the uniform continuity of $X$ and the fields $\xi_j$, for $j=1,2,3,$ we can find points $\{p_i^k\;|\;(i,k)\in I\}\in M(\mathcal{ J}_1)-\overline{M(\mathcal{ J}')}$ satisfying Statements 1), 2), 3), and the following property:
\begin{equation}\label{basesS-3}
\left\|\xi_j(p_i^k)-\xi_j(p_{i+1}^k)\right\|<\ep_0/6\mu\;,\quad \forall j \in \{1, 2,3\}\;, \quad \forall\, (i,k) \in I\;.
\end{equation}

Labeling $G$ as the spherical Gauss map of $X$, we can write $G(p_i^k)= \sum_{j=1}^3 \varrho_j^{i,k} \cdot \xi_j(p_i^k),$ $\varrho_j^{i,k} \in [-1,1].$
Take $a \in [0,1] - \{ \varrho_2^{i,k}\;|\; (i,k)\in I\}$, and define $e_1^{i,k}=\xi_1(p_i^k)$, $e_2^{i,k}=-\sqrt{1-a^2} \xi_2(p_i^k)+a \xi_3(p_i^k)$ and $e_3^{i,k}=a \xi_2(p_i^k)+\sqrt{1-a^2}\xi_3(p_i^k)$. Then, (\ref{basesS-1}) is a direct consequence of (\ref{basesS-3}). Moreover, note that $e_3^{i,k} \neq G(p_i^k)$, $\forall (i,k)$, and so (\ref{basesS-2}) trivially holds.
Finally, the existence of $\theta_i^k$ is straightforward.
\end{proof}
\begin{figure}[htbp]
	\begin{center}
		\includegraphics[width=0.45\textwidth]{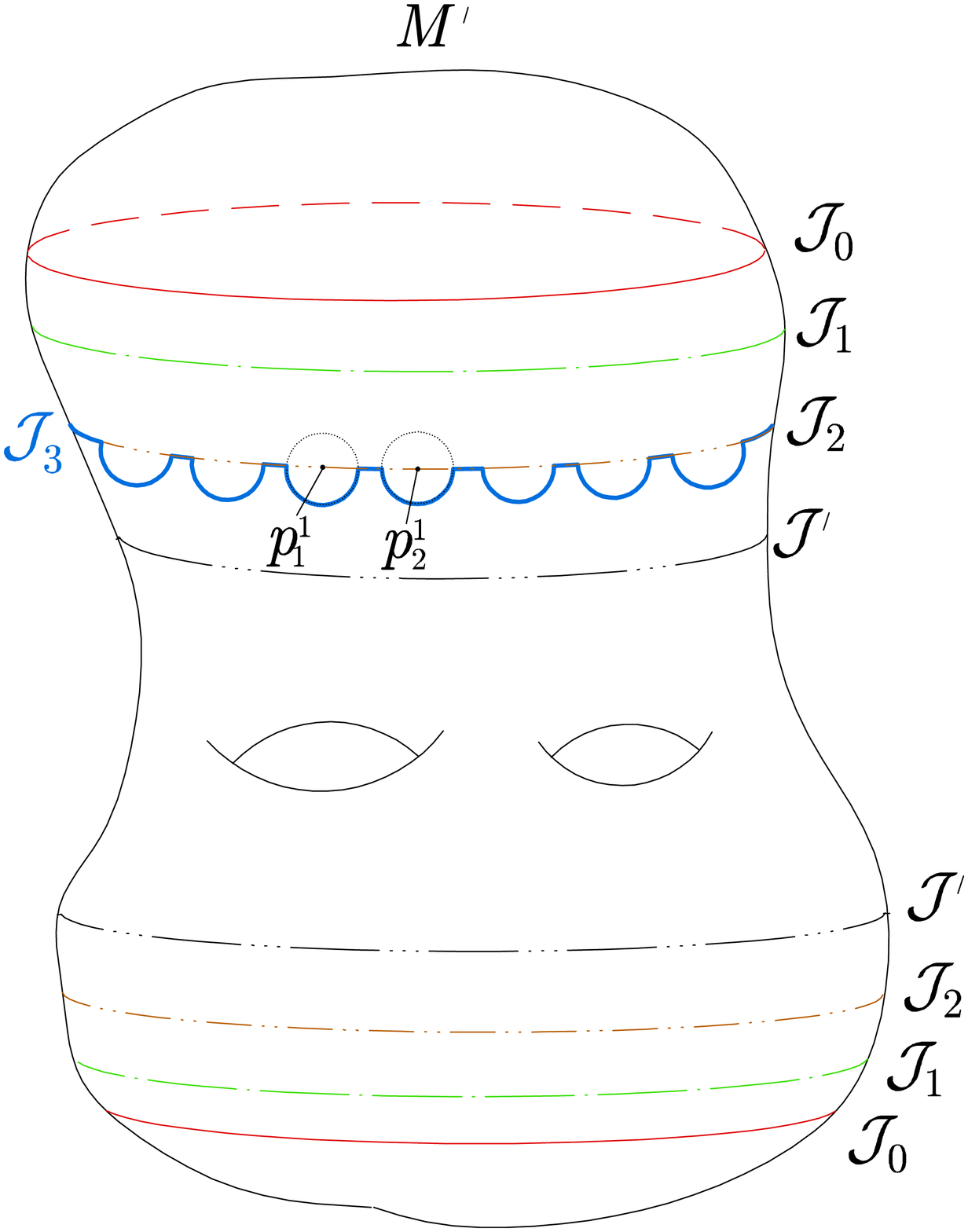}
	\end{center}
	\caption{The surface $M'$ the multicycles $\mJ_0$, $\mJ_1$, $\mJ_2$, $\mJ_3$, and $\mJ'$.}
	\label{fig:new}
\end{figure}

\begin{remark}
Notice that Properties (\ref{Bik-peque}) and (\ref{basesS-1}) are cyclic, i.e., they are true for $i=n$ labeling $p_{n+1}^k=p_1^k$, $S_{n+1}^k=S_1^k$ and $B^{n+1,k}=B^{1,k}.$
\end{remark}

Now, for any $(i,k)\in I,$ consider a holomorphic function $\zeta_{i,k}:M(\mathcal{J}_0)-\{p_i^k\}\to\C$ having a simple pole at $p_i^k.$ The existence of such functions is a consequence of the Noether `gap' Theorem (see \cite{Farkas}). Up to multiplying $\zeta_{i,k}$ by a complex constant, we can assume that the residue of $\zeta_{i,k} \cdot \omega$ at $p_i^k$ is $-1$, for all $(i,k)\in I.$

\begin{claim}\label{curvas-explicitas}
There exists $0<\de<1$ such that, for any $(i,k)\in I,$ there exist a point $q_i^k\in\partial D(p_i^k,\de)$ and a regular simple curve $\be_{i,k}:[0,1]\to \overline{D(p_i^k,\de)}$ satisfying
\begin{enumerate}
\item $\be_{i,k}(0)=q_i^k,$ $\be_{i,k}(1)=p_i^k$ and $\be_{i,k}(]0,1[)\subset D(p_i^k,\de);$
\item $\zeta_{i,k}(\be_{i,k}(t))\cdot \om_{\be_{i,k}(t)}(\be_{i,k}'(t))\in\R^+,$ $\forall\,t\in[0,1[;$
\item ${\rm Im}(\zeta_{i,k}(\be_{i,k}(t))\cdot {\rm Im}(\theta_i^k)<0,$ $\forall\,t\in[0,1[.$
\end{enumerate}
\end{claim}
At this point we can define the following constant:
$$\de':=\max\{\text{length}_{ds}(\be_{i,k})\;|\;(i,k)\in I\}\;.$$
Notice that $\de' \geq \de$ and $\lim_{\de\to 0}\de'=0.$ 

\begin{claim}\label{Af-A}
There exists $\de>0$ small enough to satisfy Claim \ref{curvas-explicitas} and the following list of properties:
\begin{enumerate}[\rm ({A}1)]

\item There exists $\mathcal{ J}_3$ a multicycle with $M(\mathcal{ J}_3)=M(\mathcal{ J}_2)-\cup_{(i,k)\in I} D(p_i^k,\de)$ (see Fig. \ref{fig:new}) ;

\item $\overline{D(p_i^k,\de)\cup D(p_{i+1}^k,\de)}\subset B^{i,k},\quad \forall (i,k)\in I$;

\item $\overline{D(p_i^k,\de)}\cap \overline{D(p_j^l,\de)}=\emptyset,\quad \forall (i,k)\neq (j,l)\in I;$

\item $\delta'\cdot\max_{\overline{D(p_i^k,\de)}}\, \{|f_{(X,S_i^k)}|\} <2\:\ep_0,\quad \forall (i,k)\in I;$

\item $\de'\cdot \max_{\overline{D(p_i^k,\delta)}} \, \{|f_{(X,S_i^k)}g^2_{(X,S_i^k)}|\} <|{\rm Im} (\theta_i^k)|\,\ep_0,$ $\forall (i,k)\in I;$

\item $3\mu \cdot \max_{p \in \overline{D(p_i^k,\delta)}} \{|f_{(X,S_i^k)}(p)-f_{(X,S_i^k)}(p_i^k)|\}<|f_{(X,S_i^k)}(p_i^k)|\,\ep_0,\quad \forall (i,k)\in I;$

\item $\delta' \cdot \max_{\overline{D(p_i^k,\delta)}}\,  \{\|\phi\|\}<\ep_0,\quad \forall (i,k)\in I$, where $\Phi=\phi\cdot\om$ is the Weierstrass representation of the immersion $X$.
\end{enumerate}
\end{claim}
Now, we label $\ell:=\text{diam}_{ds}(M(\mathcal{ J}_3))+2\:\de'+2 \, \pi \, \delta+1.$
\begin{figure}[htbp]
	\begin{center}
		\includegraphics[width=0.35\textwidth]{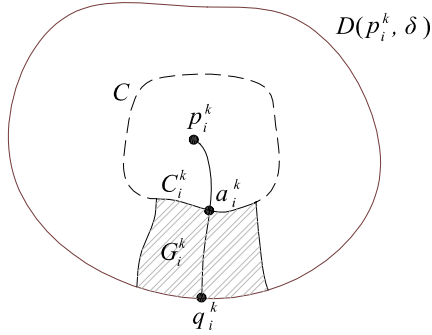}
	\end{center}
	\caption{The conformal disk $D(p_i^k,\delta)$.}
	\label{fig:disco}
\end{figure}
For each $k=1, \ldots ,\textsc{e}$, we construct a  sequence $\Psi_k=\{\Psi_{i,k}\;|\; i=1, \ldots, n\}$,  where the element $\Psi_{i,k}=\{\kappa_i^k,a_i^k,C_i^k,G_i^k,\Phi^{i,k}\}$ is composed of:
\begin{itemize}
\item $\kappa_i^k$ is a positive real number; 
\item $a_i^k$ is the first point in the (oriented) curve $\be_{i,k}$ , such that
\begin{equation}\label{ecu-a}
\frac{1}{2}\left|f_{(X,S_i^k)}(p_i^k)\right|\int_{\be(q_i^k,a_i^k)} \kappa_i^k \zeta_{i,k}\,\om=3 \mu\;,
\end{equation}
where $\mu$ was defined at the beginning of the proof of Lemma \ref{properness} and  $\be(q,p)$ denotes the oriented arc of $\be_{i,k}([0,1])$ starting at $q$ and finishing at $p.$ 
\item $C_i^k$ is a piece of a simple closed regular curve $C$ contained in $D(p_i^k,\de)$ such that $a_i^k\in C_i^k$ and each connected component of $\be_{i,k}([0,1])-\{a_i^k\}$ lies on a connected component of $\overline{D(p_i^k,\de)}-C$ (see Fig. \ref{fig:disco}.)

\item $G_i^k$ is  a closed annular sector bounded by $C_i^k$, $\partial D(p_i^k,\de)$ and the boundary of a small neighborhood of the curve $\be(q_i^k,a_i^k);$

\item $\Phi^{i,k}=\phi^{i,k}\cdot\omega$ is a Weierstrass representation defined on 
$\displaystyle  \overline{M(\mathcal{ J}_1)} - \cup_{(j,l)\leq (i,k)} U(p_j^l),$
 where $U(p_j^l)=D(p_j^l,\de)-\overline{G_j^l}$ is a small open neighborhood of $p_j^l.$
\end{itemize}

\begin{remark} \label{prime} In each family $\Psi_k$  we will adopt the convention that $\Psi_{n+1,k}= \Psi_{1,k}$. In case $k=1$, let $\Phi^{0,1}=\phi^{0,1}\omega$ be the Weierstrass representation of the immersion $X$. We denote $\Psi_{0,1}=\{\Phi^{0,1}\}$. In case $k>1$, we write $\Phi^{0,k}=\Phi^{n,k-1}$ and label $\Psi_{0,k}=\{\Phi^{0,k}\}$.

\end{remark}

\begin{claim}\label{Af-B}
We can construct the sequence in such way that satisfy 
\begin{enumerate}[\rm (B1${_i^k}$)]
\item $\de'\cdot\max_{\overline{D(p_j^l,\de)}}\{|f_{(\Phi^{i,k},S_j^l)}|\}<2\ep_0,\quad \forall\, (j,l)>(i,k);$

\item $\de'\cdot\max_{\overline{D(p_j^l,\de)}} \{|f_{(\Phi^{i,k},S_j^l)}g^2_{(\Phi^{i,k},S_j^l)}|\}<|{\rm Im}(\theta_i^k)|\,\ep_0,$ $\forall\, (j,l)>(i,k);$

\item $3 \mu \cdot \max_{p \in \overline{D(p_j^l,\de)}} \{|f_{(\Phi^{i,k},S_j^l)}(p)-f_{(X,S_j^l)}(p_j^l)|\}<|f_{(X,S_j^l)}(p_j^l)|\,\ep_0,\quad\forall\, (j,l)>(i,k);$

\item $\|{\rm Re}\int_{\alpha_p}\Phi^{i,k}\|<\ep_0,\quad \forall p\in C_i^k$, where $\alpha_p$ is a piece of $C_i^k$ connecting $a_i^k$ with $p;$

\item $\Phi_{(3,S_i^k)}^{i,k}=\Phi_{(3,S_i^k)}^{i-1,k}$, where $\Phi_{(j,S_i^k)}^{i,k}$ represents the $j$-th coordinate of the triple $\Phi^{i,k}$ in the frame $S_i^k;$

\item $\displaystyle\|\phi^{i,k}(p)-\phi^{i-1,k}(p)\|<\frac{\epsilon_0}{n\textsc{e}\ell},\quad \forall p\in K_i^k:=\overline{M(\mathcal{ J}_1)}-\left(D(p_i^k,\de) \cup \left(\bigcup_{(j,l)< (i,k)}U(p_j^l)\right)\right);$

\item $\|{\rm Re}\int_{\be(q_i^k,a_i^k)}\Phi^{i,k}-{\rm Re} \int_{\be(q_{i-1}^k,a_{i-1}^k)}\Phi^{i-1,k}\|< 15\, \epsilon_0$, (for $i=2,\ldots,n+1$);

\item For all $p\in G_i^k$ one has
\begin{multline*}
\left\|\left({\rm Re}\int_{q_i^k}^p\Phi_{(1,S_i^k)}^{i,k}\right)e_1^{i,k}+\left({\rm Re}\int_{q_i^k}^p\Phi_{(2,S_i^k)}^{i,k}\right)e_2^{i,k}-\frac{1}{2}\left|f_{(X,S_i^k)}(p_i^k)\right| \left({\rm Re}\int_{q_i^k}^p\kappa_i^k \zeta_{i,k}\,\om\right)\:e_1^{i,k}\right\|<5\ep_0\;;
\end{multline*}

\item $3 \mu +\ep_0 \geq \frac{1}{2}\big|f_{(X,S_i^k)}(p_i^k)\big| ({\rm Re}\int_{q_i^k}^p \kappa_i^k \zeta_{i,k}\,\om)\geq -\ep_0,$ for all $p\in \overline{G_i^k}.$
\end{enumerate}
\end{claim}

The above properties are true for $(i,k) \in I$, except for (B1$_i^k$), (B2$_i^k$), and (B3$_i^k$) which hold only for $(i,k) \neq (n,\textsc{e})$. Similarly, Property (B7$_i^k$) is valid only for $i=2, \ldots,n+1, $  and any $k \in \{1,\ldots,\textsc{e}\}$ (see Remark \ref{prime}.) 

\label{japon}We define each family $\Psi_k$ in a recursive way.  Before entering in the details of the recursive construction, we would like to make some remarks:

\begin{remark}To construct $\Psi_{i,1}$ starting from $\Psi_{i-1,1}$ we will use Properties
(B1$_{i-1}^1$), (B2$_{i-1}^1$) and (B3$_{i-1}^1$). In the case $i-1=0$, these properties are a consequence of (A4), (A5), and (A6), respectively.

If $k>1$ and $\Psi_{i-1,k}$ is already defined, then  in order to obtain $\Psi_{i,k}$  from $\Psi_{i-1,k}$ we will make use of Properties
(B1$_{i-1}^k$), (B2$_{i-1}^k$) and (B3$_{i-1}^k$). In the case $i-1=0$, these properties are a consequence of (B1$_{n}^{k-1}$), (B2$_{n}^{k-1}$) and (B3$_{n}^{k-1}$), respectively. \end{remark}
Assume $\Psi_{i-1,k}$ is defined satisfying Properties (B1$_{i-1}^k$),$\dots$, (B9$_{i-1}^k$). 

From item (3) in Claim \ref{curvas-explicitas} we easily obtain that:
\begin{equation} \label{eq:leonor}
\left|1+ c \: \theta_i^k \: \zeta_{i,k} (\beta_{i,k}(t)) \right| \geq |{\rm Im}(\theta_i^k)|>0, \quad \forall t \in [0,1], \quad \forall c>0.           
\end{equation}
Consider a basis of the homology of $M(\mathcal{J}_0)$, \label{base-ho}$ \mathcal{ B}=\{\aleph_1,\ldots,\aleph_\varrho\}$, so that the curves $\aleph_j$, $j=1, \ldots, \varrho,$ are contained in $M(\mathcal{J}_2)-\cup_{(j,l)\in I}D(p_j^l, \delta).$ 

Reasoning as in Claim \ref{3Runge} we obtain the existence of $\{\wect{1}{2\varrho}{\vp}\}\subset \mathcal{O}(\overline{M(\mathcal{ J}_0)})$ such that $\{\mathcal{F}(\vp_1),$ $\dots,$ $\mathcal{F}(\vp_{2\varrho})\}$ are linearly independent. Up to a suitable shrinking, we can assume
\begin{equation}\label{R2-2}
\left| \exp\left[\sum_{j=1}^{2\varrho}x_j \vp_j(p)\right]-1 \right|<\frac{|{\rm Im}(\theta_i^k)|}{2}\;,
\end{equation}
$\forall\, x=(\wect{1}{2\varrho}{x})\in\R^{2\varrho}$ with $\left\|x\right\|<1$, $\forall\, p\in\overline{M(\mathcal{ J}_0)}$.
Now, for $\T=(\lambda_0, \lambda_1,\ldots,\l_{2\varrho})\in\R^{2\varrho+1},$ we consider the map $h^{\T}:\overline{M(\mathcal{ J}_0)} \to\C$ given by
$$
h^{\T}(p)=\lambda_0 \: \theta_i^k \: \zeta_{i,k}(p)+\exp\left[\sum_{j=1}^{2\varrho}\l_j\vp_j(p) \right]\;.
$$
Observe that $h^\T \to 1$ uniformly on $\overline{M(\mathcal{ J}_0)}-D(p_i^k,\delta)$, as $\T \to 0.$ Then, there exists $1>r>0$, so that $h^\T$ has no zeroes in $\overline{M(\mathcal{ J}_0)}-D(p_i^k,\delta)$, for all $\Theta \in B(0,r)=\{x \in  \R^{2 \varrho+1}\; | \; \|x\|<r\}.$

Label $g^{\Theta}=g_{(\Phi^{i-1,k},S_i^k)}/h^{\T}$ and $\Phi_{3}^{\T}=\Phi_{(3,S_i^k)}^{i-1,k}$. For the associate Weierstrass representation, $\Phi^\T$, we define the period function $\mathcal{ P}:\R^{2\varrho+1}\to\R^{2\varrho}$ given by 
$$
\mathcal{ P}(\T)=\left( {\rm Re}\left[\int_{\aleph_j}\Phi_1^{\T}\right]_{j=1,\dots,\varrho}\, ,\, {\rm Re}\left[ \int_{\aleph_j}\Phi_2^{\T}\right]_{j=1,\dots,\varrho}\right)\;.
$$
Notice that $\mathcal{P}$ is a mapping of class  $C^1$ and $\mathcal{P}(0,\ldots,0)=0.$ Then, applying the Implicit Function Theorem, as in the proof of Lemma \ref{AdaptaRunge-2}, we get the existence of a positive constant $\kappa>0$ and a curve $L:]-\kappa,\kappa[ \to \R^{2 \varrho}$, such that $(\lambda_0,L(\lambda_0)) \in B(0,r)$ and $\mathcal{P}(\lambda_0,L(\lambda_0))=0,$ for all $\lambda_0$ in $]-\kappa,\kappa[.$ Since $\Phi^{(\lambda_0,L(\lambda_0))}\to\Phi^{i-1,k}$, uniformly on $K_i^k,$ as $\lambda_0 \to 0$, then we can find $\kappa_i^k \in ]0,\kappa[$ so that the Weierstrass data:
\begin{equation}\label{datosW}
g_{(\Phi^{i,k},S_i^k)}:=g^{(\kappa_i^k ,L(\kappa_i^k))}\;,\quad \Phi_{(3,S_i^k)}^{i,k}:=\Phi_3^{(\kappa_i^k ,L(\kappa_i^k))}\;,
\end{equation}
satisfy Properties (B1$_i^k$), (B2$_i^k$), (B3$_i^k$) and (B6$_i^k$). Furthermore, Property (B5$_i^k$) trivially follows from the definition of $\Phi^{i,k}$.

For the sake of simplicity we will write  $h_{i,k}$ instead of $h^{(\kappa_i^k ,L(\kappa_i^k))}$.  We would like to point out that the immersion $X_i^k:\overline{M(\mathcal{ J}_2)}-\left(D(p_i^k,\delta) \cup\left(\cup_{(j,l)< (i,k)}U(p_j^l)  \right) \right)\to\R^3$ with Weierstrass representation $\Phi^{i,k}$, in the orthogonal frame $S_i^k$, is well-defined. To obtain the remainder properties we have to work a little further. 

To check Property (B8$_i^k$) we write $a+{\rm i}b\equiv a e^{i,k}_1+b e^{i,k}_2$. Given $p\in \be(q_i^k,a_i^k)$, from the definition of $\be_{i,k}$ we get
$$ {\rm Re}\int_{q_i^k}^p\kappa_i^k \zeta_{i,k}\,\om=\int_{\be(q_i^k,p)}\kappa_i^k \zeta_{i,k}\,\om\in\R^+\;.
$$
Hence, using first (\ref{zetas}) and then (\ref{ecu-a}), one obtains:
$$ \frac{1}{2}\left| \left( \left| f_{(X,S_i^k)}(p_i^k) \right|\right. - 
\left.\,\overline{\t_i^k\, f_{(X,S_i^k)}(p_i^k)} \right) \int_{\be(q_i^k,p)}\kappa_i^k \zeta_{i,k}\,\om \right|<
\frac{1}{2}\,\frac{\ep_0}{3\mu} \left| f_{(X,S_i^k)}(p_i^k)\right|\int_{\be(q_i^k,p)}\kappa_i^k \zeta_{i,k}\,\om\leq \ep_0\;. $$

Therefore, we have
$$
\left|{\rm Re}\int_{\be(q_i^k,p)}\Phi_{(1,S_i^k)}^{i,k}+{\rm i} {\rm Re}\int_{\be(q_i^k,p)}\Phi_{(2,S_i^k)}^{i,k}-
\frac{1}{2}\left|f_{(X,S_i^k)}(p_i^k)\right|\int_{\be(q_i^k,p)}\kappa_i^k \zeta_{i,k}\,\om \right|<
$$
\begin{equation}\label{(**)}
\left| {\rm Re}\int_{\be(q_i^k,p)}\Phi_{(1,S_i^k)}^{i,k}+{\rm i}{\rm Re}\int_{\be(q_i^k,p)}\Phi_{(2,S_i^k)}^{i,k}-\right. 
\left.  \frac{1}{2}\,\overline{\theta_i^k f_{(X,S_i^k)}(p_i^k)}\int_{\be(q_i^k,p)}\kappa_i^k \zeta_{i,k}\,\om \right|+ \ep_0
\end{equation}
Taking into account the definition of $h_{i,k}$ and  \eqref{R2-2}, we can write $h_{i,k}=v_{i,k}+\t_i^k\kappa_i^k \zeta_{i,k}+1,$  where $v_{i,k}$ is a holomorphic function with $|v_{i,k}|< \frac{|{\rm Im}(\theta_i^k)|}{2}$. Moreover, ${\rm Re}\,\Phi_1+{\rm i}{\rm Re}\,\Phi_2=\frac{1}{2}(\overline{\eta}- g^2 \eta)$. Then, expression (\ref{(**)}) can be bounded by
\begin{multline*}
\frac{1}{2}\left| \int_{\be(q_i^k,p)}\overline{f_{(\Phi^{i-1,k},S_i^k)}\t_i^k\kappa_i^k \zeta_{i,k}\,\om}\,+ \int_{\be(q_i^k,p)}\overline{f_{(\Phi^{i-1,k},S_i^k)}(v_{i,k}+1)\,\om}\,-\right. \\
\left.\int_{\be(q_i^k,p)}f_{(\Phi^{i-1,k},S_i^k)}g^2_{(\Phi^{i-1,k},S_i^k)}\frac{\om}{h_{i,k}}\,- \overline{\theta_i^k f_{(X,S_i^k)}(p_i^k)} \int_{\be(q_i^k,p)}\kappa_i^k \zeta_{i,k}\,\om \right|+ \ep_0 \leq
\end{multline*}
$$
\frac{1}{2}\left|\int_{\be(q_i^k,p)}\overline{(f_{(\Phi^{i-1,k},S_i^k)}-f_{(X,S_i^k)}(p_i^k))\t_i^k\kappa_i^k \zeta_{i,k}\,\om} \;\right|+\frac{1}{2}\left|\int_{\be(q_i^k,p)}\overline{f_{(\Phi^{i-1,k},S_i^k)}(v_{i,k}+1)\,\om}\right|+
$$
$$
\frac{1}{2}\left|\int_{\be(q_i^k,p)}f_{(\Phi^{i-1,k},S_i^k)}g^2_{(\Phi^{i-1,k},S_i^k)}\frac{\om}{h_{i,k}}\right|+ \ep_0<\ep_0+\ep_0\left(1+\frac{|{\rm Im}(\theta_i^k)|}{2}\right)+\ep_0+\ep_0<5\ep_0\;,
$$
where in the second inequality we have used (\ref{ecu-a}), (B3$_{i-1}^k$), (B1$_{i-1}^k$), (B2$_{i-1}^k$), \eqref{eq:leonor} and  \eqref{R2-2}. Thus, we have proved that Property (B8$_i^k$) holds for all $p\in \be(q_i^k,a_i^k)$. Hence, if $C_i^k$ and $G_i^k$ are chosen close enough to $a_i^k$ and $\be(q_i^k,a_i^k)$, respectively, we obtain Properties (B4$_i^k$), (B8$_i^k$) and (B9$_i^k$).
\\

Finally, we are checking (B7$_i^k$). In order to do this, we write
\begin{multline*}
\left\|{\rm Re}\,\int_{\be(q_i^k,a_i^k)}\Phi^{i,k}-{\rm Re}\,\int_{\be(q_{i-1}^k,a_{i-1}^k)}\Phi^{i-1,k}\right\|\leq  \\
  \left\|\sum_{j=1}^2 \left[ \left({\rm Re}\,\int_{\be(q_i^k,a_i^k)}\Phi_{(j,S_i^k)}^{i,k}\right)e_j^{i,k}-\left({\rm Re}\,\int_{\be(q_{i-1}^k,a_{i-1}^k)}\Phi_{(j,S_{i-1}^k)}^{i-1,k}\right)e_j^{i-1,k} \right] \right\| +\\
  \left| \left({\rm Re}\,\int_{\be(q_i^k,a_i^k)}\Phi_{(3,S_i^k)}^{i,k}\right)e_3^{i,k}-\left({\rm Re}\,\int_{\be(q_{i-1}^k,a_{i-1}^k)}\Phi_{(3,S_{i-1}^k)}^{i-1,k}\right)e_3^{i-1,k}  \right|\;,
\end{multline*}
and we separately bound each addend. Using (B8$_i^k$), (B8$_{i-1}^k$), (\ref{ecu-a}) and (\ref{basesS-1}), we obtain
\begin{gather*}
 \left\|\sum_{j=1}^2 \left[ \left({\rm Re}\,\int_{\be(q_i^k,a_i^k)}\Phi_{(j,S_i^k)}^{i,k}\right)e_j^{i,k}-\left({\rm Re}\,\int_{\be(q_{i-1}^k,a_{i-1}^k)}\Phi_{(j,S_{i-1}^k)}^{i-1,k}\right)e_j^{i-1,k} \right] \right\|< \\
  \frac{1}{2}\left\|\left( \left|f_{(X,S_i^k)}(p_i^k)\right|\int_{\be(q_i^k,a_i^k)}\kappa_i^k \zeta_{i,k}\,\om \right) e^{i,k}_1 -  \left( \left|f_{(X,S_{i-1}^k)}(p_{i-1}^k)\right|\int_{\be(q_{i-1}^k,a_{i-1}^k)}\kappa_{i-1}^k \zeta_{i-1,k}\,\om \right) e^{i-1,k}_1 \right\|+10\ep_0 =\\
\left\| 3\mu\, e^{i,k}_1-3\mu\, e^{i-1,k}_1\right\|+10\ep_0 <\ep_0+10\ep_0=11\ep_0\;.
\end{gather*}

To bound the second addend we use (B5$_i^k$) and (B5$_{i-1}^k$) to obtain
\begin{gather} \label{eq:policia}
\left\|\left({\rm Re}\,\int_{\be(q_i^k,a_i^k)}\Phi_{(3,S_i^k)}^{i,k}\right)e_3^{i,k}- \left({\rm Re}\,\int_{\be(q_{i-1}^k,a_{i-1}^k)}\Phi_{(3,S_{i-1}^k)}^{i-1,k}\right)e_3^{i-1,k}\right\|\leq \left|{\rm Re}\,\int_{\be(q_i^k,a_i^k)}\Phi_{(3,S_i^k)}^{i-1,k}\right|+\\ \nonumber
\left|{\rm Re}\,\int_{\be(q_{i-1}^k,a_{i-1}^k)}\Phi_{(3,S_{i-1}^k)}^{i-2,k}\right| \leq \de' (\text{max}_{\overline{D(p_i^k,\de)}}\{\|\phi^{i-1,k}\|\}+\text{max}_{\overline{D(p_{i-1}^k,\de)}}\{\|\phi^{i-2,k}\|\})< 2\left(\frac{\de'\ep_0}{\ell}+\ep_0\right)<4\ep_0\;,
\end{gather}
where in the second to last inequality we have used (B6$_j^l$), $(j,l)<(i,k),$ and (A7). Therefore, Property (B7$_i^k$) holds, and so we have constructed the required sequence $\{\Psi_{i,k}\;|\;(i,k)\in I\}.$


\subsubsection{Preparing the second deformation}
Note that the Weierstrass representations $\Phi^{i,k}$ have simple poles and zeros in $M(\mathcal{J}_1)$. Our next job is to describe a domain $\mathcal{U}$ in $M(\mathcal{J}_1)$ where the above Weierstrass representations determine minimal immersions.

We can consider $\delta''>\delta$ such that $\overline{D(p_i^k,\de'')\cup D(p_{i+1}^k,\de'')}\subset B^{i,k},$ $ \forall (i,k)\in I$, and  $\overline{D(p_i^k,\de'')}\cap \overline{D(p_j^l,\de'')}=\emptyset,$ $\forall (i,k)\neq (j,l)\in I.$

Let $\alpha_{i,k}\subset D(p_i^k,\de'')- \overline{D(p_i^k,\delta)}$ be a simple curve connecting $\partial D(p_i^k,\de'') \cap \intc M(\mathcal{J}_2)$ with $q_i^k$ and finally let $N_i^k$ be a small open neighborhood of $\alpha_{i,k} \cup \beta(q_i^k,a_i^k)$ in  $\overline{G_i^k \cup (D(p_i^k,\de'')- D(p_i^k,\delta))}$. The domain $\mathcal{U}$ is defined as
$$\mathcal{U}=\left( M(\mathcal{ J}_2)-\bigcup_{(i,k)\in I} D(p_i^k,\de'') \right)\cup \left( \bigcup_{(i,k)\in I} N_i^k \right)$$

If $\de''$, $\alpha_{i,k}$ and $N_i^k$ are suitably chosen, then we can  guarantee:
\begin{claim}\label{Af-C}
The domain $\mathcal{U}$ satisfies the following properties:
\begin{enumerate}[\rm ({C}1)]
\item There exists $\mathcal{ J}_\mathcal{U}$ a multicycle with $\mathcal{U}=M(\mathcal{ J}_\mathcal{U}).$ From now on, we write $M(\mathcal{ J}_\mathcal{U})$ instead of $\mathcal{U};$
\item $\be(q_i^k,a_i^k)\subset\overline{M(\mathcal{ J}_\mathcal{U})}$ and $\mathcal{J}'<\mathcal{ J}_\mathcal{U};$
\item $\text{\rm diam}_{ds}(M(\mathcal{ J}_\mathcal{U}))<\ell;$
\item $\overline{M(\mathcal{ J}_\mathcal{U})}\cap \overline{D(p_i^k,\de)}\subset \overline{G_i^k},$ $\forall\, (i,k)\in I;$
\item The homology group of $M(\mathcal{ J}_\mathcal{U})$ is the same as $M(\mathcal{ J}_0)$ and it is generated by the basis $\mathcal{B}$ described in page \pageref{base-ho}.
\end{enumerate}
\end{claim}
\noindent At this point, it is clear that we are able to find a multicycle, $\mathcal{ J}_{4}$, with $\mathcal{ J}_\mathcal{U}<\mathcal{ J}_{4}$ and satisfying (C3) and (C5), where the immersions $X_i^k:M(\mathcal{ J}_{4})\to\R^3$ given by $X_i^k(p)={\rm Re}\,\int_{p_0}^p \Phi^{i,k}$ are still well-defined, for $(i,k)\in I.$
\begin{claim}\label{Af-D}
For $(i,k)\in I,$ we have
\begin{enumerate}[\rm (D1${_i^k}$)]
\item $\|X_i^k(p)-X_{i-1}^k(p)\|<\frac{\ep_0}{n\textsc{e}},$ $\forall\, p\in M(\mathcal{ J}_{4})-D(p_i^k,\de)$;
\item $(X_i^k)_{(3,S_i^k)}=(X_{i-1}^k)_{(3,S_i^k)}$;
\item $\|X_n^\textsc{e}(a_i^k)-X_n^\textsc{e}(a_{i+1}^k)\|<20\ep_0$;
\item $X_n^\textsc{e}(a_i^k)\in\R^3- E_{2\mu}$.
\end{enumerate}
\end{claim}
\begin{proof}
In order to get (D1$_i^k$) we use (B6$_i^k$) and (C3) as follows:
\begin{multline*}
\| X_i^k(p)-X_{i-1}^k(p)\| =\left\| {\rm Re}\, \int_{p_0}^p( \phi^{i,k}-\phi^{i-1,k} )\om \right\|\leq 
\int_{p_0}^p\left|\phi^{i,k}-\phi^{i-1,k}\right|\,\|\om\| \leq \frac{\ep_0}{n\textsc{e}\ell} \int_{p_0}^p\|\om\|<\frac{\ep_0}{n\textsc{e}}\;.
\end{multline*}

Now, (B5$_i^k$) immediately implies (D2$_i^k$). To check (D3$_i^k$) we apply (D1$_j^l$), $(j,l)\in I$, (B7$_{i+1}^k$) and (\ref{Bik-peque}) to obtain
$$
\|X_n^\textsc{e}(a_i^k)-X_n^\textsc{e}(a_{i+1}^k)\| \leq \|X_n^\textsc{e}(a_i^k)-X_i^k(a_i^k)\| + \| X_n^\textsc{e}(a_{i+1}^k)-X_{i+1}^k(a_{i+1}^k) \| +
$$
$$
\| X_{i+1}^k(q_{i+1}^k)-X_i^k(q_i^k)\| + \|(X_i^k(a_i^k)-X_i^k(q_i^k))-(X_{i+1}^k(a_{i+1}^k)-X_{i+1}^k(q_{i+1}^k))\|<
$$
$$
4\ep_0+\|X(q_{i+1}^k)-X(q_i^k)\| +\left\| \mbox{Re} \left( \int_{\be(q_{i+1}^k,a_{i+1}^k)}\Phi^{i+1,k}-\int_{\be(a_i^k,q_i^k)}\Phi^{i,k}\right) \right\| <4\ep_0+\ep_0+15\ep_0=20\ep_0\;.
$$

Finally, we will prove (D4$_i^k$). Using (D1$_j^l$), $(j,l)>(i,k),$ one gets
$$
\| X_n^\textsc{e}(a_i^k)-X(p_i^k)-3\mu\,\mathcal{N}_E(X(p_i^k))\|\leq\|X_n^\textsc{e}(a_i^k)-X_i^k(a_i^k)\|+ \|X_i^k(a_i^k)-X_i^k(q_i^k)-3\mu\,\mathcal{N}_E(X(p_i^k))\|+
$$
$$
\|X_i^k(q_i^k)-X(p_i^k)\|<
\ep_0+\|(X_i^k(a_i^k)-X_i^k(q_i^k))_{(*,S_i^k)}-3\mu e^{i,k}_1\|+|(X_i^k(a_i^k)-X_i^k(q_i^k))_{(3,S_i^k)}|+
$$
$$
\|X_i^k(q_i^k)-X(q_i^k)\|+\|X(q_i^k)-X(p_i^k)\|<
\ep_0+5\ep_0+2 \:\ep_0+\ep_0+\ep_0=10\:\ep_0\;,
$$
where in the last inequality we have used (D1$_i^k$), (B8$_i^k$), (\ref{Bik-peque}), (\ref{ecu-a}) and \eqref{eq:policia}. As $X(p_i^k)+3\,\mu \, \mathcal{N}_E(X(p_i^k))\in\R^3- E_{3\mu}$, then (D4$_i^k$) holds for a small enough $\ep_0$.
\end{proof}


\subsubsection{The second deformation}

 For any $(i,k)\in I,$ let $T_i^k=\{w_1^{i,k},w_2^{i,k},w_3^{i,k}\}$ be a new  orthonormal basis such that
\begin{equation}\label{basesT}
w_3^{i,k}=\mathcal{N}_E(X_n^\textsc{e}(a_i^k))\;.
\end{equation}
 Consider also $Q_i^k$ the connected component of the set $\overline{\mathcal{ J}_\mathcal{U}-(C_i^k\cup C_{i+1}^k) \,}$ that does not cut $C_j^l$, $\forall(j,l)\in I- \{(i,k),(i+1,k)\}$. Note that $\{Q_i^k\;|\;(i,k)\in I\}$ satisfy:
\begin{equation}\label{Q1}
\overline{Q_i^k}\cap \overline{Q_j^l}=\emptyset\;, \; \text{for all } (i,k)\neq (j,l) \quad \text{and} \quad Q_i^k\subset B^{i,k}\;, \; \text{for all $(i,k) \in I$};
\end{equation}
\begin{equation}\label{Q2}
Q_i^k\cap \overline{D(p_j^l,\de)}=\emptyset\;,\quad (j,l)\notin\{(i,k),(i+1,k)\}\;
\end{equation}
and, up to a small perturbation,
\begin{equation}\label{Q3}
f_{(X_n^\textsc{e},T_i^k)}(p)\neq 0\;,\quad \forall\, p\in Q_i^k\;.
\end{equation}
 Now, let $\widehat C_i^k$ be an open set containing $C_i^k$ and sufficiently small to fulfill
\begin{equation}\label{Q4}
\|X_n^\textsc{e}(p)-X_n^\textsc{e}(a_i^k)\|<3\epsilon_0\;,\quad \forall\, p\in \widehat C_i^k\cap\overline{M(\mathcal{ J}_\mathcal{U})}\;.
\end{equation}
Notice that the above choice is possible due to Properties (D1$_j^l$), $(j,l)>(i,k),$ and (B4$_i^k$).
We also define, for any $\xi>0,$ $Q_i^k(\xi)=\{p\in M(\mathcal{J}_0) \; |\; \dist_{(M(\mathcal{J}_0),ds)}(p,Q_i^k)\leq\xi\}.$
\begin{figure}
	\begin{center}
		\includegraphics[width=0.70\textwidth]{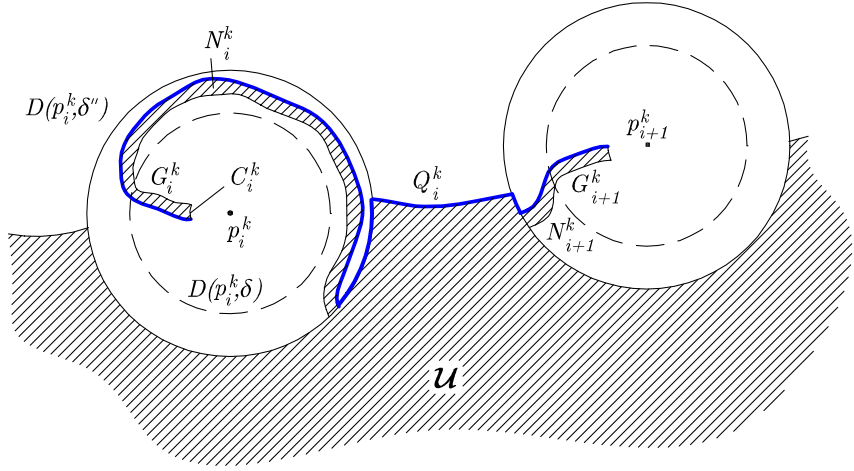}
	\end{center}
	\caption{The curves $Q_i^k.$}
	\label{fig:QQ}
\end{figure}

\begin{claim}\label{Af-E}
There exists $\xi>0$ small enough so that:
\begin{enumerate}[\rm ({E}1)]
\item $Q_i^k(\xi) \subset M(\mathcal{ J}_{4});$

\item $Q_i^k(\xi)\cap Q_j^l(\xi)=\emptyset$, for $(i,k)\not=(j,l);$

\item $Q_i^k(\xi)\cap\overline{D(p_j^l,\de)}=\emptyset$, for $(j,l)\not\in\{(i,k),(i+1,k)\};$

\item $Q_i^k(\xi)\subset B^{i,k};$

\item $|\int_{[x_0,x_1]}\om|=\xi/2,$ $\forall\, x_0\in Q_i^k,\forall\, x_1 \in \partial D(x_0,\xi/2),\forall\, (i,k)\in I$ where  $[x_0,x_1]$ represents the length minimizing arc joining $x_0$ and $x_1$ (recall that $ds^2=\|\omega\|^2$.)

\item Given $p \in Q_i^k$, we have
$|f_{(X_n^\textsc{e},T_i^k)}(p)-f_{(X_n^\textsc{e},T_i^k)}(q)|<\varepsilon_1,$ $ \forall q\in D(p,\xi/2),$  where $\varepsilon_1=\frac{1}{4}\min_{Q_i^k}\{|f_{(X_n^\textsc{e},T_i^k)}|\};$

\item $\text{\rm diam}_{ds}\left(\overline{M(\mathcal{ J}_\mathcal{U})-Q_i^k(\xi)}\right)<\ell$.
\end{enumerate}
\end{claim}
Observe that Properties {(E2)}, {(E3)}, {(E4)} and {(E7)} are consequence of (\ref{Q1}), (\ref{Q2}), and {(C3)}. Furthermore, (E5) holds as the developing map, $\mathsf{f}$, is a local isometry (see Remark \ref{re:developing}.) The other ones are straightforward.
\\

We are now ready to construct a sequence $\{\Lambda_{i,k}\;|\;(i,k)\in I\}$ where the element $\Lambda_{i,k}=\{Y_i^k,\tau_i^k,\nu_i^k\}$ is composed of:
\begin{itemize}
\item $Y_i^k:\overline{M(\mathcal{ J}_{4})}\to\R^3$ is a conformal minimal immersion. We also label $Y_0^1=X_n^\textsc{e}$ and $Y_0^k=Y_n^{k-1}, $ $k \geq 2;$
\item $\{(\tau_i^k,\nu_i^k)\in\R^+\times\R^+\;|\;(i,k)\in I\}.$
\end{itemize}

\begin{claim}\label{Af-F}
We can construct the sequence $\{\Lambda_{i,k}\;|\;(i,k)\in I\}$ satisfying the following list of properties:
\begin{enumerate}[\rm (F1${_i^k}$)]
\item $(Y_i^k)_{(3,T_i^k)}=(Y_{i-1}^k)_{(3,T_i^k)}$;
\item $\|Y_i^k(p)-Y_{i-1}^k(p)\|<\frac{\ep_0}{n\textsc{e}}$, $\forall p\in\overline{M(\mathcal{ J}_\mathcal{U})-Q_i^k(\xi)}$;
\item $|f_{(Y_i^k,T_j^l)}(p)-f_{(Y_{i-1}^k,T_j^l)}(p)|<\frac{\varepsilon_1}{n\textsc{e}}$, $\forall p\in \overline{M(\mathcal{ J}_\mathcal{U})-Q_i^k(\xi)}$, for $(j,l)>(i,k);$
\item $\left(\frac{1}{\tau_i^k}+\frac{\nu_i^k}{\tau_i^k(\tau_i^k-\nu_i^k)}\right)\text{\rm max}_{Q_i^k(\xi)}\{|f_{(Y_{i-1}^k,T_i^k)}g^2_{(Y_{i-1}^k,T_i^k)}|\}+\nu_i^k\text{\rm max}_{Q_i^k(\xi)}\{|f_{(Y_{i-1}^k,T_i^k)}|\} <\frac{2}{\xi}$;
\item $\frac{1}{2} \left( \frac{\tau_i^k \, \xi}{4} \; \text{\rm min}_{Q_i^k}\{|f_{(Y_0^1,T_i^k)}|\}-1\right)>\text{\rm diam}_{\R^3} (E')+1;$
\end{enumerate}
\end{claim}

Assume we have constructed $Y_0^1,Y_1^1,\ldots,Y_{i-1}^k$. Then we use Lemma \ref{AdaptaRunge-2} to get a holomorphic function without zeros $l_{i,k}:\overline{M(\mathcal{ J}_\mathcal{U})}\to\C$ such that
\begin{itemize}
\item $|l_{i,k}(p)-\tau_i^k|<\nu_i^k$, $\forall\, p\in Q_i^k(\xi/2)$;
\item $|l_{i,k}(p)-1|<\nu_i^k$, $\forall\, p\in\overline{M(\mathcal{ J}_\mathcal{U})-Q_i^k(\xi)}$;
\item The minimal immersion $Y_i^k$ with Weierstrass data given by
\begin{equation}\label{datosW2}
g_{(Y_i^k,T_i^k)}=\frac{g_{(Y_{i-1}^k,T_i^k)}}{l_{i,k}}\quad\text{and}\quad \Phi^{i,k}_{(3,T_i^k)}=\Phi^{i-1,k}_{(3,T_i^k)},
\end{equation}
is well-defined.
\end{itemize}
Then, we define the immersion $Y_i^k$ as $Y_i^k(p)={\rm Re}\,\int_{p_0}^p\Phi^{i,k}$, where the Weierstrass data $\Phi^{i,k}$, in the orthogonal frame $T_i^k$, are determined by the López-Ros transformation (\ref{datosW2}).  Notice that $ \phi_{(Y_i^k,T_j^l)} \stackrel{\nu_i^k \to 0}{\longrightarrow} \phi_{(Y_{i-1}^k,T_j^l)}$ uniformly on $\overline{M(\mathcal{ J}_\mathcal{U})-Q^{k,\xi}_i}$. At this point, if $\nu_i^k$ is small enough and $\tau_i^k$ is sufficiently large, then to check that $Y_i^k$ satisfies Properties (F1$_i^k$),...,(F5$_i^k$) is a straightforward computation, and so we have constructed the sequence $\{\Lambda_{i,k}\;|\;(i,k)\in I\}.$ Note that \eqref{Q3} is used in the proof of (F5$_i^k$).


\subsubsection{The immersion $Y$ solving Lemma \ref{properness}}

Consider the minimal immersion $Y:\overline{M(\mathcal{ J}_\mathcal{U})}\rightarrow \R^3$ given by $Y=Y_n^\textsc{e}$. We are going to check that $Y$ satisfies all the statements of Lemma \ref{properness}. 
\\

{\bf Item (L2.b):}
Items 2 and 3 in Claim \ref{puntos-p} and Properties {(E4)} and {(A2)} imply that
$
\overline{M(\mathcal{J}')}\subset M(\mathcal{ J}_\mathcal{U})- \big(\cup_{(i,k)\in I} D(p_i^k,\de)\big)\cup\big(\cup_{(i,k)\in I} Q_i^k(\xi)\big)\;.
$
So, we can successively apply (D1$_i^k$) and (F2$_i^k$), $(i,k)\in I,$ to obtain $\forall p\in\overline{M(\mathcal{J}')}$
\begin{equation}\label{(19)}
\|Y(p)-X(p)\|\leq\|Y_n^\textsc{e}(p)-Y_0^1(p)\|+\|X_n^\textsc{e}(p)-X(p)\|< 2\ep_0<b_1\;,
\end{equation}
where the last inequality occurs if $\ep_0$ is small enough.
\\


{\bf Items (L2.a) and (L2.c):} As a previous step we will prove the following claim:

\begin{claim}\label{2proper}
Every connected curve $\g$ in $M(\mathcal{ J}_\mathcal{U})$ connecting $M(\mathcal{J}')$ with $\mathcal{ J}_\mathcal{U}$ contains a point $p'\in\g$ such that $Y(p')\in\R^3- E'$.
\end{claim}

\begin{proof}
Let $\g\subset\overline{M(\mathcal{ J}_\mathcal{U})}$ be a connected curve with $\g(0)\in M(\mathcal{J}')$ and $\g(1)=x_0\in\mathcal{ J}_\mathcal{U}$.
\\

\noindent{\bf Case i)} Assume $x_0\in \widehat C_i^k\cap Q_i^k(\xi)$.
Using Properties (E2), (F2$_j^l$) for $(j,l)\neq (i,k)$, (F1$_i^k$) and Inequality (\ref{Q4}), we infer
\begin{multline}\label{(*)}
|(Y_n^\textsc{e}(x_0)-X_n^\textsc{e}(a_i^k))_{(3,T_i^k)}| \leq \|Y_n^\textsc{e}(x_0)-Y_i^k(x_0)\| +|(Y_i^k(x_0)-Y_{i-1}^k(x_0))_{(3,T_i^k)}|+ \\ \|Y_{i-1}^k(x_0)-Y_0^1(x_0)\|+\|X_n^\textsc{e}(x_0)-X_n^\textsc{e}(a_i^k)\|< \ep_0+\ep_0+3\ep_0=5\ep_0\;.
\end{multline}
If we write $T$ as the tangent plane to $\partial E$ at the point $\mathcal{ P}_E(X_n^\textsc{e}(a_i^k))$, then we know that $\dist_{\R^3}(p,\partial E)\geq \dist_{\R^3}(p,T)$ for any  $p$ in the halfspace determined by $T$ that does not contain $\partial E$. If $\ep_0$ is small enough, (D4$_i^k$), \eqref{basesT}, and \eqref{(*)} guarantee that $Y_n^\textsc{e}(x_0)$ belongs to the above halfspace, and moreover we have
\begin{multline}\label{(20)}
\dist_{\R^3}(Y_n^\textsc{e}(x_0), \partial E) \geq \dist_{\R^3} (Y_n^\textsc{e}(x_0),T)=(Y_n^\textsc{e}(x_0)-\mathcal{ P}_E(X_n^\textsc{e}(a_i^k)))_{(3,T_i^k)}>\\ (X_n^\textsc{e}(a_i^k)-\mathcal{ P}_E(X_n^\textsc{e}(a_i^k)))_{(3,T_i^k)}-5 \epsilon_0 >  2 \mu-5 \epsilon_0>\mu\;.
\end{multline}
From the definition of $\mu$ we conclude $Y_n^\textsc{e}(x_0)\in \R^3 - E'$.
\\

\noindent{\bf Case ii)} Assume $x_0\in \widehat C_i^k\cap Q_{i-1}^k(\xi)$.
Reasoning as in the above case and using Property (D3$_{i-1}^k$), we obtain
$$
|(Y_n^\textsc{e}(x_0)-X_n^\textsc{e}(a_{i-1}^k))_{(3,T_{i-1}^k)}|\leq  |(Y_n^\textsc{e}(x_0)-X_n^\textsc{e}(a_{i}^k))_{(3,T_{i-1}^k)}| + \|X_n^\textsc{e}(a_i^k)-X_n^\textsc{e}(a_{i-1}^k)\|<25\ep_0\;.
$$
Now, following the arguments of (\ref{(20)}), we conclude $Y(x_0)\in\R^3- E'$.
\\

\noindent{\bf Case iii)} Assume $x_0\in \widehat C_i^k- \cup_{(j,l)\in I} Q_j^l(\xi)$. Taking into account (F2$_j^l$),  for $(j,l)\in I,$ and (\ref{Q4}), one has
$$
\|Y_n^\textsc{e}(x_0)-X_n^\textsc{e}(a_i^k)\|\leq \|Y_n^\textsc{e}(x_0)-Y_0^1(x_0)\|+\|X_n^\textsc{e}(x_0)-X_n^\textsc{e}(a_i^k)\|<4\epsilon_0\;,
$$
and then we can finish as in the preceding cases.
\\

\noindent{\bf Case iv)} Finally, suppose that $x_0 \in Q_i^k - \cup_{(j,l)\in I} \widehat{C}_j^l$. For the sake of simplicity, we will write $f^{i-1,k}$ and  $g^{i-1,k}$ instead of $f_{(Y_{i-1}^k,T_i^k)}$ and $g_{(Y_{i-1}^k,T_i^k)}$, respectively, and $a+ {\rm i}\, b$ instead of $a w_1^{i,k}+b w_2^{i,k}$. Hence, for $x_1\in\g\cap\partial D(x_0,\xi/2)$, taking into account (F2$_j^l$), for $(j,l)>(i,k)$, and the definition of $Y_i^k$  one has
$$
\|Y_n^\textsc{e}(x_0)-Y_n^\textsc{e}(x_1)\|>  \|Y_i^k(x_0)-Y_i^k(x_1)\|-2\epsilon_0\geq
\|(Y_i^k(x_0)-Y_i^k(x_1))_{(*,T_i^k)}\|-2\epsilon_0=$$
$$\frac{1}{2}\left|\int_{[x_1,x_0]} \overline{f^{i-1,k}l_{i,k}\,\om}-\int_{[x_1,x_0]}\frac{f^{i-1,k}(g^{i-1,k})^2}{l_{i,k}}\,\om\right|-2\epsilon_0 \geq
\frac{1}{2}\left|\tau_i^k\int_{[x_1,x_0]} \overline{f^{i-1,k}\,\om}\right|-$$ $$
\frac{1}{2}\left|\frac{1}{\tau_i^k}\int_{[x_1,x_0]} f^{i-1,k}(g^{i-1,k})^2\,\om\right|-
\frac{1}{2}\left|\int_{[x_1,x_0]} \overline{f^{i-1,k}(l_{i,k}-\tau_i^k)\,\om}\right|-
$$
$$\frac{1}{2}\left|\int_{[x_1,x_0]} f^{i-1,k}(g^{i-1,k})^2\left(\frac{1}{l_{i,k}}- \frac{1}{\tau_i^k}\right)\,\om\right|-2\epsilon_0\geq
$$
using the definition of $l_{i,k}$ and (E5), we obtain
$$
\geq\frac{\tau_i^k}{2}\left|\int_{[x_1,x_0]} \overline{f^{i-1,k}\,\om}\right|-\frac{\xi}{4}\left(\frac{1}{\tau_i^k}\text{max}_{Q_i^k(\xi)}\{| f^{i-1,k}(g^{i-1,k})^2|\}+\nu_i^k\text{max}_{Q_i^k(\xi)}\{|f^{i-1,k}|\}+\right.
$$
$$
\left.
\frac{\nu_i^k}{\tau_i^k(\tau_i^k-\nu_i^k)}\text{max}_{Q_i^k(\xi)}\{|f^{i-1,k}(g^{i-1,k})^2|\}\right)-2\epsilon_0\geq \frac{1}{2}\left(\tau_i^k\left|\int_{[x_1,x_0]} \overline{f^{i-1,k}\,\om}\right|-1\right)-2\epsilon_0\;,
$$
where we have used (F4$_i^k$) in the last inequality.
On the other hand, we make use of  (E5), (E6), and (F3$_j^l$), $(j,l)<(i,k),$ to deduce
$$
\left|\int_{[x_1,x_0]} f^{i-1,k\,}\om\right|\geq\left|f_{(Y_0^1,T_i^k)}(x_0)\int_{[x_1,x_0]} \om\right|-
\left|\int_{[x_1,x_0]} (f_{(Y_0^1,T_i^k)}(x_0)-f_{(Y_0^1,T_i^k)})\,\om\right|
$$
$$
-\left|\int_{[x_1,x_0]} (f_{(Y_0^1,T_i^k)}-f^{i-1,k})\,\om\right|\geq 
\tfrac{\xi}2\, (|f_{(Y_0^1,T_i^k)}(x_0)|-2 \varepsilon_1 )\geq
\tfrac{\xi}4 \, \text{min}_{Q_i^k}\{|f_{(Y_0^1,T_i^k)}|\}\;.
$$
Therefore, by using (F5$_i^k$) for $\ep_0$ small enough we have
$$
\|Y_n^\textsc{e}(x_0)-Y_n^\textsc{e}(x_1)\|> \frac{1}{2} \left( \tau_i^k \tfrac{\xi}4 \, \text{min}_{Q_i^k}\{|f_{(Y_0^1,T_i^k)}|\}-1\right)-2\epsilon_0>
\text{diam}_{\R^3} (E')+1-2\epsilon_0>\text{diam}_{\R^3} (E')\;.
$$
From the above inequality we conclude that $\g$ satisfies the claim in this last case.
It is clear that $x_0$ has to lie in one of the above cases, hence, we have proved the claim.
\end{proof}

Moreover, if $\ep_0$ is small enough, (\ref{(19)}) and the convex hull property for minimal surfaces guarantee that $Y(\overline{M(\mathcal{J}')})\subset E'$. Claim \ref{2proper} implies that we can find a multicycle $\mathcal{ J}$ satisfying (L2.a) and (L2.c).
\\

{\bf Item (L2.d):} Given $p\in \overline{M(\mathcal{ J})}- M(\mathcal{J}')$ there are five possible situations for the point $p$ (recall that $Q_i^k(\xi) \cap \overline{D(p_j^l,\de)}=\emptyset$, $(j,l)\notin\{(i,k),(i+1,k)\}).$
\\

\noindent{\bf Case I)} Suppose $p\not\in (\cup_{(i,k)\in I} D(p_i^k,\de))\cup(\cup_{(i,k)\in I} Q_i^k(\xi))$. In this case we can use  Properties (D1$_i^k$), (F2$_i^k$), $(i,k)\in I$ to conclude that: 
$$\|Y_n^\textsc{e}(p)-X(p)\|\leq \|Y_n^\textsc{e}(p)-Y_0^1(p)\| + \|X_n^\textsc{e}(p)-X(p)\|<\ep_0+\ep_0=2\ep_0<2b_2\;.
$$
As usual, we have assumed that $\ep_0$ is small enough.

The above fact jointly with Hypothesis (\ref{(3)}) of Lemma \ref{properness} give us that $Y(p) \not\in E_{-2 b_2}$.
\\

\noindent{\bf Case II)} Suppose $p\in D(p_i^k,\de)- \cup_{(j,l)\in I} Q_j^l(\xi)$, for an $(i,k)\in I$.
In this case, one has
\begin{multline*}
\left< Y_{n}^\textsc{e}(p)-X(p_i^k),e_1^{i,k}\right> = \left< Y_n^\textsc{e}(p)-Y_0^1(p), e^{i,k}_1 \right> + \left< X_{n}^\textsc{e}(p)-X_i^k(p),e^{i,k}_1 \right> +\\
\left< X_i^k(p)-X_i^k(q_i^k),e^{i,k}_1 \right>+ \left< X_i^k(q_i^k)-X(q_i^k),e^{i,k}_1\right>  + \left< X(q_i^k)-X(p_i^k),e^{i,k}_1\right> >
\end{multline*}
using (D1$_j^l$), $\forall\, (j,l)\neq (i,k)$, (F2$_j^l$), $\forall\,(j,l)\in I$, and (\ref{Bik-peque}),
$$
> \left< X_i^k(p)-X_i^k(q_i^k),e^{i,k}_1\right>-4\ep_0> 
 \frac{1}{2}|f_{(X,S_i^k)}(p_i^k)| \left({\rm Re}\int_p^{q_i^k}\kappa_i^k \zeta_{i,k}\,\om\right)-9\epsilon_0\geq -10 \epsilon_0> -b_2\;,
$$
where we have used (B8$_i^k$) and (B9$_i^k$). Recall that $e_1^{i,k}=\mathcal{N}_E(X(p_i^k))$. Therefore, again
as a consequence of Hypothesis (\ref{(3)}), we infer $Y_{n}^\textsc{e}(p)\not\in E_{-b_2}$. In particular  $Y_{n}^\textsc{e}(p)\not\in E_{-2 b_2}$.
\\

\noindent{\bf Case III)} Assume $p\in D(p_i^k,\de)\cap Q_i^k(\xi)$, for some $(i,k)\in I.$ This case is slightly more complicated.

As a previous step we need to get an upper bound for $\| w_3^{i,k}-e_1^{i,k} \|$. Remember that when we checked (D4$_i^k$), we obtained
$\|X_n^\textsc{e}(a_i^k)- (3\mu e_1^{i,k}+X(p_i^k) ) \| \leq 11\epsilon_0.$
Therefore,
\begin{multline}\label{(21)}
\| w_3^{i,k}-e_1^{i,k} \|=\| \mathcal{N}_E(X_n^\textsc{e}(a_i^k))-\mathcal{N}_E(X(p_i^k) )\|=\\
\| \mathcal{N}_E(X_n^\textsc{e}(a_i^k))-\mathcal{N}_E(3\mu e_1^{i,k}+X(p_i^k) )\| \leq
M \|X_n^\textsc{e}(a_i^k)- (3\mu e_1^{i,k}+X(p_i^k) ) \| \leq 11M \epsilon_0\;,
\end{multline}
where $M$ represents the maximum of $\|d\mathcal{N}_E\|$ in $\R^3- E$. Note that $M$ does not depend on $\ep_0$.
On the other hand, using (F1$_i^k$) and (\ref{(21)}), we find
\begin{equation}\label{(22)}
\left|\left< Y_i^k(p)-Y_{i-1}^k(p),e_1^{i,k}\right>\right|=  \left|\left< Y_i^k(p)-Y_{i-1}^k(p),e_1^{i,k}-w_3^{i,k}\right>\right|\leq 11M \ep_0 (\|Y_i^k(p)\|+\|Y_{i-1}^k(p)\|)\;.
\end{equation}

Now, making use of (F2$_j^l$), $(j,l)\neq (i,k)$, (\ref{(22)}) and (D1$_j^l$), $(j,l)>(i,k),$ one obtains
\begin{multline}\label{caso3}
\left< Y_n^\textsc{e}(p)-X(p_i^k),e_1^{i,k}\right>\geq \left< Y_i^k(p)-X(p_i^k),e_1^{i,k}\right>-\ep_0\geq
\left< Y_{i-1}^k(p)-X(p_i^k),e_1^{i,k}\right>- \\
11M \ep_0 (\|Y_i^k(p)\|+\|Y_{i-1}^k(p)\|)-\ep_0\geq \left< X_i^k(p)-X(p_i^k),e_1^{i,k}\right>-11M \ep_0 (\|Y_i^k(p)\|+\|Y_{i-1}^k(p)\|)-3\ep_0\;.
\end{multline}

At this point, we can argue as in the previous case to conclude
\begin{equation}\label{(23)}
\left< Y_n^\textsc{e}(p)-X(p_i^k),e_1^{i,k}\right> > -b_2-11M \ep_0 (\|Y_i^k(p)\|+\|Y_{i-1}^k(p)\|)-3\ep_0\;.
\end{equation}
Observe that Item (L2.c), the convex hull property and (F2$_j^l$), $(j,l) > (i,k),$ guarantee that $Y_i^k(p)\in E'_{\ep_0}$. Furthermore, notice that
\begin{multline}\label{(24)}
\|Y_{i-1}^k(p)-X(q_i^k)\|\leq \|Y_{i-1}^k(p)-X_i^k(p)\|+\| (X_i^k(p)-X(q_i^k))_{(*,S_i^k)}\|+ \\
|(X_i^k(p)-X(q_i^k))_{(3,S_i^k)}|<2\ep_0+\ep_0+5\ep_0+3\mu+ \ep_0+2 \ep_0=3\mu +11\ep_0\;,
\end{multline}
where we have used (F2$_j^l$), $(j,l)<(i,k),$ and (D1$_j^l$), $(j,l)>(i,k),$ to get a bound of the first addend; (B8$_i^k$), (B9$_i^k$) and (D1$_j^l$), $(j,l) \leq (i,k),$ to get a bound of the second addend; and (D2$_i^k$), (D1$_j^l$), $(j,l)< (i,k),$ and (\ref{Bik-peque}) to get a bound of the third one. Then $\|Y_i^k(p)\|$ and $\|Y_{i-1}^k(p)\|$ are bounded in terms of $\ep_0$. So, we infer from (\ref{(23)}) that  $Y(p)\not\in E_{-2b_2}$, if $\ep_0$ is small enough.
\\

\noindent{\bf Case IV)} Suppose  $p\in D(p_{i+1}^k,\de)\cap Q_i^k(\xi)$. Reasoning as in the preceding case, now we can deduce from \eqref{basesS-1} $\| e_1^{i+1,k}-w_3^{i,k} \| \leq \| e_1^{i+1,k}-e_1^{i,k}\|+\|e_1^{i,k}-w_3^{i,k}\|<\frac{\ep_0}{3 \mu}+11M\ep_0$ and obtain
$$
\left|\left< Y_i^k(p)-Y_{i-1}^k(p),e_1^{i+1,k}\right>\right|=\left|\left< Y_i^k(p)-Y_{i-1}^k(p),e_1^{i+1,k}-w_3^{i,k}\right>\right|\leq 
(11M\ep_0+\tfrac{\ep_0}{3 \mu} )(\|Y_i^k(p)\|+\|Y_{i-1}^k(p)\|)\;.
$$
Using these inequalities as in the former case, we deduce $Y(p)\not\in E_{-2b_2}$.
\\

\noindent{\bf Case V)} Finally, assume $p\in Q_i^k(\xi)-\cup_{(j,l)\in I} D(p_j^l,\de)$.
Reasoning as in inequality \eqref{caso3}, we have
$$
\left< Y_n^\textsc{e}(p)-X(p_i^k),e_1^{i,k}\right> > \left< X(p)-X(p_i^k),e_1^{i,k}\right>-11 M \ep_0 (\|Y_i^k(p)\|+\|Y_{i-1}^k(p)\|)-3\ep_0\;,
$$
and using now (\ref{Bik-peque}), we obtain for a sufficiently small $\ep_0$,
$$\left< Y_n^\textsc{e}(p)-X(p_i^k),e_1^{i,k}\right> >-11 M \ep_0 (\|Y_i^k(p)\|+\|Y_{i-1}^k(p)\|)-4\ep_0\geq-2b_2\;.$$

This concludes the proof of Item (L2.d) and completes the proof of Lemma \ref{properness}.

\begin{remark}\label{k1>0}
If $E$ is strictly convex, then the above proof also gives that
$$
\|Y(p)-X(p)\|<\mathcal{M}(b_2,E,E'):=\sqrt{ \frac{2\left(\de^H(E,E')+2 \, b_2\right)}{\kappa_1(\partial E)}+\de^H(E,E')^2}\;,\quad \forall p\in \overline{M(\mathcal{ J})}-M(\mathcal{J}')\;,
$$
where $\de^H$ means the Hausdorff distance.
\end{remark}




\section{Completeness Lemmas}

This is the moment of employing the Runge type result proved in Section \ref{sec:runge} as well as López-Ros deformation in order to perturb a given minimal surface with finite topology about its boundary. In this way, we are able of increasing the intrinsic diameter of the surface, but preserving the extrinsic one. 
The proofs of the lemmas bellow are inspired in a new technique introduced by Nadirashvili and the last author in \cite{plateau}.

In order to state the next lemma, we shall denote $M=M' - \cup_{i=1}^\textsc{e} \D_i$, where $\D_i$, $i=1, \ldots, \textsc{e}$, are conformal disks in the compact surface $M'$. As in the previous section, $\omega$ will represent a holomorphic 1-form without zeros in $M$ and $ds^2=\|\omega\|^2.$ 
For any $i\in\{1,\ldots,\textsc{e}\},$ let $\Sigma_i$ be an analytic cycle around $\D_i$ and $\be_i:\Sigma_i\to\Gamma_i\subset\R^3$ an analytic Jordan curve. Given $\mathcal{T}(\Sigma_i)$ a tubular neighborhood of $\Sigma_i$ in $(M,ds^2)$, we denote by ${\tt P}_i:\mathcal{T}(\Sigma_i) \to \Sigma_i$ the natural projection. In this setting we have:

\begin{lemma}\label{nadira}
Consider $\mathcal{ J}=\{\g_1,\dots,\g_\textsc{e}\}$ a multicycle  on $M$, $X:\overline{M(\mathcal{ J})}\to\R^3$ a conformal minimal immersion, $p_0$ a point in $M(\mathcal{ J})$, and $r>0$,  such that: 
\begin{enumerate}[\rm (1)]
\item $X(p_0)=0$
\item $\gamma_i \subset \mathcal{T}(\Sigma_i),$ for $i=1, \ldots, \textsc{e};$
\item $\|X(p)-\be_i({\tt P}_i(p))\|<r,$  for all $p\in \g_i$ and for all $i=1,\ldots,\textsc{e}\;.$
\end{enumerate}
Then, for any $s>0,$ and any $\ep>0$ so that $p_0\in M(\mathcal{ J}^\ep)$,  there exist $\widetilde{\mathcal{ J}}=\{\widetilde{\g}_1,\dots,\widetilde{\g}_\textsc{e}\}$ a multicycle and a conformal minimal immersion $\widetilde{X}: \overline{M(\widetilde{\mathcal{ J}})} \to \R^3$, with $\widetilde X(p_0)=0$, and satisfying:
\begin{enumerate}[\rm (L\ref{nadira}.a)]
\item $\widetilde \gamma_i \subset \mathcal{T}(\Sigma_i),$ for $i=1, \ldots, \textsc{e};$
\item $\mathcal{ J}^\ep<\widetilde{\mathcal{ J}}<\mathcal{ J}$;

\item $s<\dist_{(\overline{M(\widetilde{\mathcal{ J}})},\widetilde{X})}(p,\widetilde{\mathcal{ J}}),$ $\forall p\in \mathcal{ J}^\ep$;

\item $\|\widetilde{X}(p)-\be_i({\tt P}_i(p))\|<R=\sqrt{4s^2+r^2}+\ep,$ $\forall p\in \widetilde{\g}_i,$ $\forall i=1,\ldots,\textsc{e}.$ 
\end{enumerate}
\end{lemma}


\begin{lemma}\label{main}
Let $\mathcal{ J}=\{\g_1,\ldots,\g_\textsc{e}\}$ be a multicycle, $X:\overline{M(\mathcal{ J})}\to\R^3$ a conformal minimal immersion, and $p_0$ a point in $M(\mathcal{ J})$ such that $X(p_0)=0.$

Then, for any $\l>0$ and for any $\mu>0$ so that $p_0 \in M(\mathcal{ J}^\mu)$, there exists a multicycle $\widehat{\mathcal{ J}}=\{\widehat{\g}_1,\ldots,\widehat{\g}_\textsc{e}\}$ and a conformal minimal immersion $\widehat{X}:\overline{M(\widehat{\mathcal{ J}})}\to \R^3$, with $\widehat X(p_0)=0$, and  satisfying:
\begin{enumerate}[\rm (L\ref{main}.a)]
\item $\mathcal{ J}^\mu<\widehat{\mathcal{ J}}<\mathcal{ J};$

\item $\dist_{(\overline{M(\widehat{\mathcal{ J}})},\widehat{X})}(p,\widehat{\mathcal{ J}})> \l,$ $\forall p\in \mathcal{ J}^\mu;$

\item $\|X-\widehat{X}\|<\mu,$ in $M(\widehat{\mathcal{ J}}).$
\end{enumerate}
\end{lemma}


\subsection{Proof of Lemma \ref{nadira}}

As analytic Jordan curves are dense in the set of piecewise regular Jordan curves, we can assume (without lost of generality) that the multicycle $\mathcal{J}$ is analytic. Let $\zeta_0\in]0,\ep[$ be small enough so that $\gamma_i^{\zeta_0} \subset \mathcal{T}(\Sigma_i),$ for $i=1, \ldots,\textsc{e}$. Consider $N\in\N$  such that $2/N<\zeta_0, $ and:
\begin{equation} \label{eq:(a)} \left\{ \begin{array}{r}
\|X(p)-\be_i({\tt P}_i(p))\|<r, \text{ for all $p$ in the connected component of} \\ \text{ $\overline{M(\mathcal{ J})}-M(\mathcal{ J}^{2/N})$ around $\D_i$, $\forall i=1,\ldots,\textsc{e}.$} \end{array} \right.
\end{equation}
\begin{remark}
Throughout the proof of the lemma a set of real positive constants  depending on $X$, $\mathcal{ J}$, $r$,  $\ep$, and $s$ will appear. The symbol `$\cte$' will denote these different constants. It is important to note that the choice of these constants does not depend on $N$.
\end{remark}

For the sake of simplicity, we will consider again an order relation in the set $I\equiv \{1,\ldots,2N\}\times\{1,\ldots,\textsc{e}\}.$ We say $(j,l)>(i,k)$ if one of the two following situations occurs: $l=k$ and $j>i$ or $l>k.$

For each $k=1,\ldots,\textsc{e},$ let $\{v_{1,k}, \ldots ,v_{2N,k}\}$ be a set of points in the curve $\g_k$ that divide $\g_k$ into $2N$ equal parts (i.e., curves with the same length). Following the normal projection, we can transfer the above partition to the curve $\g_k^{2/N}$: $\{v_{1,k}^\prime, \ldots, v_{2N,k}'\}$. We define the following sets:

\begin{itemize}
\item $L_{i,k}=[v_{i,k}\,,\,v_{i,k}']$, $\forall\, (i,k)\in I.$ Recall that $[v_{i,k}\,,\,v_{i,k}']$ represents the minimizing geodesic in $(M(\mJ),ds^2)$ joining $v_{i,k}$ and $v_{i,k}'$;

\item $\mathcal{ G}_{j,k}=\g_k^{j/N^3}$, $\forall\, j=0,\ldots, 2N^2$ (recall that $\g_k^{j/N^3}$ means the parallel curve to $\g_k$, in $M(\mJ),$ such that the distance between them is $j/N^3$);

\item $\mathcal{ A}_k=\bigcup_{j=0}^{N^2-1}\overline{\Int{\mathcal{ G}_{2j+1,k}}- \Int{\mathcal{ G}_{2j,k}}}$ and $\widetilde{\mathcal{ A}}_k=\bigcup_{j=1}^{N^2}\overline{\Int{\mathcal{ G}_{2j,k}}- \Int{\mathcal{ G}_{2j-1,k}}};$

\item $\mathcal{ R}_k= \bigcup_{j=0}^{2N^2} \mathcal{ G}_{j,k};$

\item $\mathcal{ B}_k= \bigcup_{j=1}^N L_{2j,k}$ and $\widetilde{\mathcal{ B}}_k= \bigcup_{j=0}^{N-1} L_{2j+1,k};$

\item $\mathcal{ L}_k=\mathcal{B}_k \cap \mathcal{A}_k$, $\mathcal{ \widetilde{L}}_k=\widetilde{\mathcal{B}}_k \cap \widetilde{\mathcal{A}}_k$, and $H_k=\mathcal{R}_k \cup \mathcal{ L}_k \cup \mathcal{ \widetilde{L}}_k$;

\item $\Om_{N,k}=\{ p \in \Int{(\mathcal{ G}_{2N^2,k})} - \Int{(\mathcal{ G}_{0,k})}\;|\; \dist_{(M,ds)}(p,H_k) \geq \frac{1}{4N^3}\}$;

\item $\Om_N=\bigcup_{k=1}^\textsc{e} \Om_{N,k};$

\item $\om_i^k$ is the union of the curve $L_{i,k}$ and those connected components of $\Om_{N,k}$ that have nonempty intersection with $L_{i,k}$ for $(i,k)\in I;$

\item $\varpi_i^k = \{ p \in M \; | \; \dist_{(M,ds)}(p,\om_i^k)< \de(N) \}$, where  $\de(N)>0$ is chosen in such a way that the sets $\overline{\varpi}_i^k,$ $(i,k)\in I,$ are pairwise disjoint.
\end{itemize}
\begin{figure}[htbp]
	\begin{center}
		\includegraphics[width=.85\textwidth]{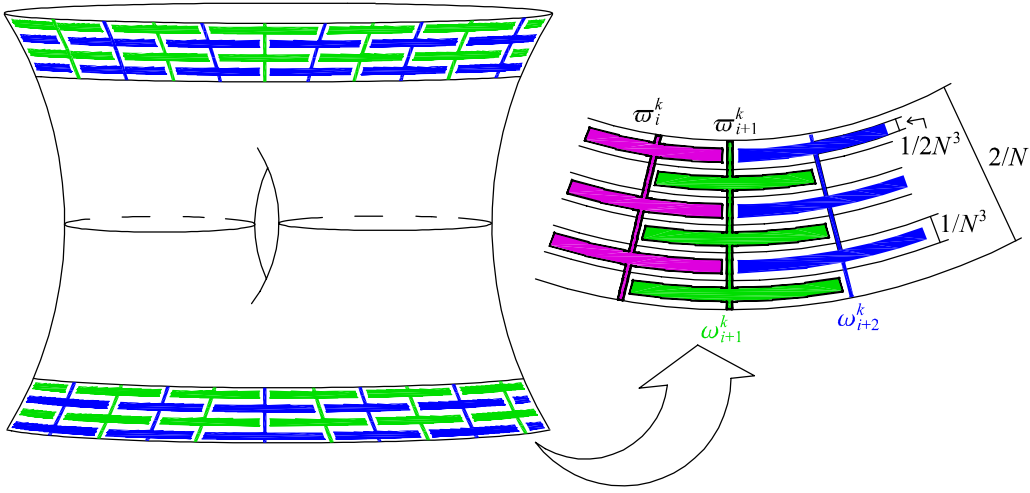}
	\end{center}
	\caption{The labyrinth around the boundary of $M(\mJ).$}
	\label{fig:laberinto}
\end{figure}

\begin{claim}\label{1nadira}
If $N$ is large enough, for any $(i,k)\in I,$  one has
\begin{enumerate}[\rm  (1)]
\item\label{ome-peque} ${\rm diam}_{(M,ds)}(\varpi_i^k)<\frac{\cte}{N}$;

\item\label{met-gran} If $\l^2\cdot ds^2$ is a conformal metric on $\overline{M(\mathcal{ J})}$ that satisfies
$$\l\geq
\begin{cases}
c & \text{ in } M(\mathcal{ J})\\
c\; N^4 & \text{ in } \Om_N\;,
\end{cases}$$
for $c>0$, and if $\a$ is a curve in $\overline{M(\mathcal{ J})}$ connecting $\gamma_k^{\zeta_0}$ and $\gamma_k$, for some $k \in \{1, \ldots,\text{e}\}$, then we have $\longui_{(M,\l\cdot ds)}(\a)\geq \cte c \; N$.
\end{enumerate}
\end{claim}
\begin{proof} The proof of item (\ref{ome-peque}) in the above claim  is straightforward. In order to prove item (\ref{met-gran}),  we denote $\alpha_j$ as the piece of $\alpha$ connecting $\gamma_k^{j/N}$ and $\gamma_k^{(j+1)/N}$, for $j=0, \ldots, N^2-1$. Then either the length of $\alpha_j$ (in $(M,ds^2)$) is greater than $\frac{\cte}{N}$ or the  length of $\alpha_j \cap \Omega_{N,k}$ is greater than $\frac{1}{2 N^3}.$ To see the former assertion, the reader only have to consider that this fact is true for curves in $\C$ and take into account that the developing map of $\omega$ is a local isometry (see Remark \ref{re:developing}.) These facts and our assumption about $\l$ give us item (\ref{met-gran}). \end{proof}

At this point, for a sufficiently large $N$, we construct a sequence of conformal minimal immersions (with boundary) defined on $\overline{M(\mathcal{ J})}$,
$\{F_i^k\;|\;(i,k)\in I\},$
by using López-Ros transformations with parameters given by Lemma \ref{AdaptaRunge-2}. We consider $F_0^1=X$ and denote $F_0^k=F_{2N}^{k-1},$ $\forall k=2,\ldots,\textsc{e}.$
\begin{claim}\label{Af-a}
These immersions will be constructed to satisfy
\begin{enumerate}[\rm ({b}1$_{i}^{k}$)]
\item $F_i^k(p)= {\rm Re}\, \left( \int_{p_0}^p \Phi^{i,k}\right)$, where $\Phi^{i,k}=\phi^{i,k}\,\om;$

\item $\|\phi^{i,k}(p)-\phi^{i-1,k}(p) \| \leq 1/N^2,$ for all $p \in \overline{M(\mathcal{ J})}- \varpi_i^k$;

\item $\|\phi^{i,k}(p) \| \geq N^{7/2},$ for all $p \in \om_i^k$;

\item $\| \phi^{i,k}(p) \| \geq \frac{\cte}{\sqrt{N}},$ for all $p \in \varpi_i^k$;

\item $\dist_{\esf^2}(G_i^k(p), G_{i-1}^k(p))<\frac{1}{N^2},$ for all $p \in \overline{M(\mathcal{ J})}-\varpi_i^k$, where $\dist_{\esf^2}$ is the intrinsic distance in $\esf^2$ and $G_i^k$ represents the Gauss map of the immersion $F_i^k$;

\item There exists an orthonormal  basis of $\R^3,$ $S_i^k=\{ e_1, e_2, e_3 \}$ such that
\begin{enumerate}[\rm ({b6}.1$_{i}^{k}$)]
\item For any $p \in \overline{\varpi_i^k}$ with $\|X(p)-\be_k({\tt P}_k(p))\|\geq 1/\sqrt{N}$, we have $\| (X(p)-\be_k({\tt P}_k(p)))_{(*,S_i^k)} \| <\frac{\cte}{\sqrt{N}}$;
\item $(F_{i}^k(p))_{(3,S_i^k)}=(F_{i-1}^k(p))_{(3,S_i^k)},$ for all $p \in \overline{M(\mathcal{ J})}$;
\end{enumerate}

\item $\| F_i^k(p)-F_{i-1}^k(p) \| \leq \frac{\cte}{N^2}$, $\forall\, p \in \overline{M(\mathcal{ J})} -\varpi_i^k.$
\end{enumerate}
\end{claim}
\begin{proof}
The sequence $\{F_i^k\;|\;(i,k)\in I\}$ is constructed in a recursive way. The order we will follow in this recursive construction is similar to the procedure explained in page \pageref{japon} for the family $\Psi_{i,k}$. When $i-1=0$ we adopt the convention that $F_0^k:=F_{2N}^{k-1}$, if $k>1$, and $F_0^1:=X$. The same occurs for the Weiertrass representations. 

Suppose that we have $\{F_j^l\;|\; (j,l)<(i,k)\}$ satisfying Items (b1$_j^l)$, $ \ldots,$ (b7$_j^l$). First we need to check the following assertions.

\begin{claim}\label{Af-b}
For a  large enough $N$, the following statements hold:
\begin{enumerate}[\rm ({c}1)]
\item $\|\phi^{i-1,k}\| \leq \cte$ in $\overline{M(\mathcal{ J})} - \cup_{(j,l)<(i,k)} \varpi_i^k$;

\item $\|\phi^{i-1,k}\| \geq \cte$ in $\overline{M(\mathcal{ J})} - \cup_{(j,l)<(i,k)} \varpi_i^k$;

\item The diameter in $\R^3$ of $F_{i-1}^k(\varpi_i^k)$ is less than $\frac{1}{\sqrt{N}}$;

\item The diameter in $\esf^2$ of $G_{i-1}^k(\varpi_i^k)$ is less than $\frac{1}{\sqrt{N}}$. In particular $G_{i-1}^k(\varpi_i^k) \subset \hbox{Cone}\left(g , \frac{1}{\sqrt{N}}\right)$, for some $g \in G_{i-1}^k(\varpi_i^k)$, where  $\hbox{Cone}(x,\theta):=\{y \in \R^3 \; | \; \angle(x,y)<\theta \};$

\item There exists an orthogonal frame $S_i^k=\{e_1, e_2, e_3 \}$ in $\R^3$, satisfying
\begin{enumerate}[\rm ({c5.}1)]
\item $\angle (e_3, X(p)-\be_k({\tt P}_k(p))) \leq \frac{\cte}{\sqrt{N}},$ for all $p \in \varpi_i^k$ with $\|X(p)-\be_k({\tt P}_k(p))\|\geq 1/\sqrt{N};$
\item $\angle (\pm e_3,G_{i-1}^k(p)) \geq \frac{\cte}{\sqrt{N}},$ for all $p \in \varpi_i^k$.
\end{enumerate}
\end{enumerate}
\end{claim}

To deduce (c1), we write $\|\phi^{i-1,k}\|\leq \sum_{(j,l)<(i,k)}\|\phi^{j,l}-\phi^{j-1,l}\|+\|\phi^{0,1}\|\leq 2\textsc{e}/N+\|\phi^{0,1}\|\leq\cte$, where we have used (b2$_j^l$), $(j,l)<(i,k).$ Using the same property and taking $N$ large enough, we have $\|\phi^{i-1,k}\|\geq \|\phi^{0,1}\|-\sum_{(j,l)<(i,k)}\|\phi^{j,l}-\phi^{j-1,l}\|\geq \|\phi^{0,1}\|-2\textsc{e}/N\geq\cte$, so we have obtained Property (c2). To check (c3), consider $p,p'\in\varpi_i^k$, then
$$
\|F_{i-1}^k(p)-F_{i-1}^k(p')\|=\Big\|\int_p^{p'}\phi^{i-1,k}\,\om \Big\| \leq \int_p^{p'}\|\phi^{i-1,k}\|\,\|\om\| \leq 
 \cte \cdot \text{diam}_{(M,ds)}(\varpi_i^k)<\frac{\cte}{N}<\frac{1}{\sqrt{N}}\;,
$$
where we have used (c1), Claim \ref{1nadira}.\ref{ome-peque} and we have taken $N$ large enough. Now, observe that using Claim \ref{1nadira}.\ref{ome-peque} we obtain $\text{diam}_{\esf^2}(G_0^1(\varpi_i^k))<\sup\{\|(dG_0^1)_p\|\;|\; p\in \varpi_i^k\}\,\text{diam}_{(M,ds)}(\varpi_i^k)<\frac{\cte}{N},$ therefore, (b5$_j^l$), $(j,l)<(i,k),$ guarantee (c4). Finally, in order to prove (c5), consider $C=\hbox{Cone}\left(g,\frac{2}{\sqrt N}\right),$ where $g$ is given by Property (c4), and
$$\mathcal{ N}=\left\{\left. \frac{X(p)-\be_k({\tt P}_k(p))}{\|X(p)-\be_k({\tt P}_k(p))\|}\;\right|\; p\in \varpi_i^k\text{ and }\|X(p)-\be_k({\tt P}_k(p))\|\geq 1/\sqrt{N}\right\}\;.
$$
To obtain (c5.2) it suffices to take $e_3$ in  $\esf^2 - H$, where
$H=(-C) \cup C.$
On the other hand, in order to satisfy (c5.1), the vector $e_3$ must be chosen as follows:
\begin{itemize}
\item If $(\esf^2- H)\cap \mathcal{ N}\not=\emptyset$, then we take $e_3$ in that set;
\item If $(\esf^2- H)\cap \mathcal{ N}=\emptyset$, then we take $e_3\in \esf^2 - H$ satisfying $\angle(e_3,q')<\frac{2}{\sqrt N}$ for some $q'\in \mathcal{ N}.$
\end{itemize}
It is straightforward to check that this choice of $e_3$ guarantees (c5). \end{proof}

At this point we are able to construct the element $F_i^k$. Let $(g^{i-1,k},\Phi_3^{i-1,k})$ be the Weierstrass data of $F_{i-1}^k$ in the frame $S_i^k$. Applying Lemma \ref{AdaptaRunge-2}, we can construct a family of holomorphic functions $h_\a:\overline{M(\mathcal{ J})}\to\C^*$ satisfying
\begin{itemize}
\item $|h_\a-\a|<1/\a$, in $\om_i^k$;
\item $|h_\a-1|<1/\a$, in $\overline{M(\mathcal{ J})}- \varpi_i^k;$
\item The minimal immersion $F_i^k(p)={\rm Re}\, \int_{p_0}^p \Phi^{i,k}$ is well-defined in $\overline{M(\mathcal{ J})}$,
\end{itemize}
where $\a >0.$ Using $h_\a$ as a López-Ros parameter, we define the Weierstrass data of $F_i^k$ in the coordinate system $S_i^k$ as $g^{i,k}=g^{i-1,k}/h_\a$ and $\Phi^{i,k}_3=\Phi_3^{i-1,k}$. 
Taking into account the fact that $h_\alpha \to 1$ (resp. $h_\alpha \to \infty$) uniformly on $\overline{M(\mathcal{ J})} -\varpi_i^k$ (resp. on $\om_i^k$), as $\alpha \to \infty$, it is clear that properties (b1$_i^k$), (b2$_i^k$), (b3$_i^k$), (b5$_i^k$), and (b7$_i^k$) hold for a large enough value of the parameter $ \a$. Moreover, (b6.2$_i^k$) trivially holds and (b6.1$_i^k$) is a immediate consequence of (c5.1). In order to prove (b4$_i^k$), observe that from (c5.2) we obtain
$$
\frac{\sin \left(\frac{\cte}{\sqrt{N}}\right)}{1+\cos\left(\frac{\cte}{\sqrt{N}}\right)} \leq |g^{i-1,k}| \leq \frac{\sin \left(\frac{\cte}{\sqrt{N}}\right)}{1-\cos\left(\frac{\cte}{\sqrt{N}}\right)} \qquad \hbox{in } \varpi_i^k\;,
$$
and so, taking (c2) into account one has (if $N$ is large enough)
$$
\| \phi^{i,k} \| \geq | \phi^{i,k}_3|=| \phi^{i-1,k}_3| \geq \sqrt{2}\| \phi^{i-1,k} \| \frac{|g^{i-1,k}|}{1+|g^{i-1,k}|^2}\geq \cte \cdot \sin\left(\tfrac{\cte}{\sqrt{N}}\right) \geq \tfrac \cte{\sqrt{N}} \qquad \hbox{in } \varpi_i^k\;.
$$

\begin{proposition}\label{Af-c}
If $N$ is large enough, then $F_{2N}^\textsc{e}$ satisfies
\begin{enumerate}[  \rm ({d}1)]
\item $2s< \dist_{(\overline{M(\mathcal{ J})},F_{2N}^\textsc{e})}(\mathcal{ J},\mathcal{ J}^{\zeta_0})$;

\item $\|F_{2N}^\textsc{e}(p)-X(p) \| \leq \frac{\cte}{N},$ $\forall\, p\in  \overline{M(\mathcal{ J})} - \cup_{(i,k)\in I} \varpi_i^k$;

\item There exists a multicycle $\widetilde{\mathcal{ J}}=\{\widetilde{\g}_1,\ldots,\widetilde{\g}_\textsc{e}\}$ satisfying

\begin{enumerate}[  \rm ({d3}.1)]
\item $\mathcal{ J}^{\zeta_0}<\widetilde{\mathcal{ J}}<\mathcal{ J}$;

\item $s <\dist_{(\overline{M(\mathcal{ J})},{F_{2N}^\textsc{e}})}(p,M(\mathcal{ J}^{\zeta_0}))<2s$, $\forall\, p \in \widetilde{\mathcal{ J}}$;

\item The curve $\widetilde \gamma_i \subset \mathcal{T}(\Sigma_i),$ for $i=1, \ldots, \textsc{e};$

\item $\|F_{2N}^\textsc{e}(p)-\be_k({\tt P}_k(p))\|< R,$ $\forall\, p\in \widetilde{\g}_k,$ $\forall k=1,\ldots,\textsc{e}.$
\end{enumerate}
\end{enumerate}
\end{proposition}

\begin{proof}
Properties (c2), (b2$_i^k$), (b3$_i^k$) and (b4$_i^k$), $(i,k)\in I,$ guarantee
$$
\|\phi^{2N,\textsc{e}}\|\geq
\begin{cases}
\frac{\cte}{\sqrt{N}} & \text{in } M(\mathcal{ J}) \\
\frac{\cte}{\sqrt{N}} N^4 & \text{in } \Om_N\;.
\end{cases}
$$
Moreover, we know $ds_{F_{2N}^\textsc{e}}^2=\frac{1}{2}\|\phi^{2N,\textsc{e}}\|^2\, ds^2$. Therefore, if $N$ is large enough,  from Claim \ref{1nadira}.\ref{met-gran} we have 
$$
\dist_{(\overline{M(\mathcal{ J})},{F_{2N}^\textsc{e}})}(\mathcal{ J},\mathcal{ J}^{\zeta_0})\geq\cte\frac{\cte}{\sqrt{N}}N=\cte\sqrt{N}>2s\;,
$$ 
which proves item (d1). Property (d2) is deduced from (b7$_i^k$), $(i,k)\in I.$

In order to construct the multicycle $\widetilde{\mathcal{ J}}$ of the statement (d3), we consider the set 
$$
\mathcal{ D}=\{p\in \overline{M(\mathcal{ J})}- M(\mathcal{ J}^{\zeta_0})\;|\; s<\dist_{\overline{(M(\mathcal{ J})},F_{2N}^\textsc{e})}(p,M(\mathcal{ J}^{\zeta_0}))<2s\}\;.
$$
From (d1), $\mathcal{ D}\neq\emptyset$ and $\mathcal{ J}$ and $\mathcal{ J}^\ep$ are contained in different connected components of $M-\mathcal{ D}.$ Therefore, we can choose a multicycle $\widetilde{\mathcal{ J}}$ on $\mathcal{ D}$ satisfying (d3.1), (d3.2) and (d3.3).

The proof of (d3.4) is more complicated. Consider $k\in\{1,\ldots,\textsc{e}\},$ $q\in \widetilde{\g}_k$ and assume that $F_{2N}^\textsc{e}(q)\neq\be_k({\tt P}_k(q)),$ otherwise we have nothing to prove. At this point, we have to distinguish two cases:
\\

\noindent{\bf Case 1.} Suppose $q\notin\cup_{(i,k)\in I}\varpi_i^k.$ Then, item (d2) gives $\|F_{2N}^\textsc{e}(q)-X(q)\|\leq\cte/N.$ Hence, taking \eqref{eq:(a)} into account and choosing $N$ large enough we obtain $\|F_{2N}^\textsc{e}(q)-\be_k({\tt P}_k(q))\|\leq r<R.$
\\

\noindent{\bf Case 2.} Suppose there exists $(i,k)\in I$ with $q\in\varpi_i^k.$ In this situation, item  (d3.2) guarantees the existence of a curve $\zeta:[0,1]\to M(\mathcal{ J})$ satisfying $\zeta(0)\in \mathcal{ J}^\ep$, $\zeta(1)=q$ and $\text{length}(\zeta,{F_{2N}^\textsc{e}})\leq 2s$. Label $\overline{t}=\sup \{ t \in [0,1] \; | \; \zeta(t) \in \partial \varpi_i^k \}$ and $\overline{q}=\zeta(\overline{t}).$
Notice that the previous supremum exists because  $\varpi_i^k \subset \overline{M(\mathcal{ J})} - M(\mathcal{ J}^\ep)$ (for a large enough $N$). Then, taking Properties (b7$_j^l$), $(j,l)>(i,k),$ into account, we obtain
\begin{multline}\label{26}
\| F_i^k(\overline{q})-F_i^k(q)\| \leq \| F_i^k(\overline{q})-F_{2N}^\textsc{e}(\overline{q}) \| + \|F_{2N}^\textsc{e}(\overline{q})-F_{2N}^\textsc{e}(q)\| +\|F_{2N}^\textsc{e}(q)-F_i^k(q)\| \leq \\
\leq \frac{\cte}{N}+\text{length}(\zeta,{F_{2N}^\textsc{e}})+\frac{\cte}{N} \leq \frac{\cte}{N}+2s\;.
\end{multline}
On the other hand, using again (b7$_j^l$), for $(j,l)>(i,k)$, one has
\begin{equation}\label{8ii}
\|F_{2N}^\textsc{e}(q)-\be_k({\tt P}_k(q))\|\leq\|F_i^k(q)-\be_k({\tt P}_k(q))\|+\frac{\cte}{N}\;.
\end{equation}
Once more, we have to discuss two different cases:

\noindent{\bf Case 2.1.} Assume $\|X(q)-\be_k({\tt P}_k(q))\|\leq 1/\sqrt{N}.$ Hence, using \eqref{26}, (b7$_j^l$), for $(j,l) \leq (i,k)$, and (c3), we get
\begin{multline*}
\|F_i^k(q)-\be_k({\tt P}_k(q))\|\leq \|F_i^k(q)-F_i^k(\overline{q})\|+\|F_i^k(\overline{q})-F_{i-1}^k(\overline{q})\|+
\|F_{i-1}^k(\overline{q})-F_{i-1}^k(q)\|+\\ \|F_{i-1}^k(q)-X(q)\|+
\|X(q)-\be_k({\tt P}_k(q))\|\leq \frac{\cte}{N}+2s+\frac{\cte}{N^2}+\frac{1}{\sqrt{N}}+\frac{\cte}{N}+\frac{1}{\sqrt{N}}< R\;,
\end{multline*}
where $N$ has to be large enough. The above inequality and (\ref{8ii}) gives (d3.3).
\\

\noindent{\bf Case 2.2.} Assume now that $\|X(q)-\be_k({\tt P}_k(q))\|>1/\sqrt{N}.$ Then we can use (b6.2$_i^k$), (b7$_j^l$), for $(j,l)<(i,k),$ and \eqref{eq:(a)} to obtain
\begin{multline}\label{9ii}
|(F_i^k(q)-\be_k({\tt P}_k(q)))_{(3,S_i^k)}|=|(F_{i-1}^k(q)-\be_k({\tt P}_k(q)))_{(3,S_i^k)}|\leq \\
|(F_{i-1}^k(q)-X(q))_{(3,S_i^k)}|+|(X(q)-\be_k({\tt P}_k(q)))_{(3,S_i^k)}|\leq\frac{\cte}{N}+r\;.
\end{multline}

On the other hand, using \eqref{26}, (b7$_j^l$), for $(j,l) \leq (i,k)$, (c3) and (b6.1$_i^k$) one has
\begin{multline}\label{10ii}
\|(F_i^k(q)-\be_k({\tt P}_k(q)))_{(*,S_i^k)}\|\leq \|(F_i^k(q)-F_i^k(\overline{q}))_{(*,S_i^k)}\|+ \|(F_i^k(\overline{q})-F_{i-1}^k(\overline{q}))_{(*,S_i^k)}\|+ \\
\|(F_{i-1}^k(\overline{q})-F_{i-1}^k(q))_{(*,S_i^k)}\|+ \|(F_{i-1}^k(q)-X(q))_{(*,S_i^k)}\|+ \|(X(q)-\be_k({\tt P}_k(q)))_{(*,S_i^k)}\|\leq \\
\frac{\cte}{N}+2s+\frac{\cte}{N^2}+\frac{1}{\sqrt{N}}+\frac{\cte}{N}+\frac{\cte}{\sqrt{N}}\leq 2s+\frac{\cte}{\sqrt{N}}\;.
\end{multline}

Therefore, making use of (\ref{9ii}) and (\ref{10ii}), we infer
$$
\|F_i^k(q)-\be_k({\tt P}_k(q))\|<\sqrt{\left(2s+\frac{\cte}{\sqrt{N}} \right)^2+\left(r+\frac{\cte}{N} \right)^2}\;.
$$
Then, using this upper bound and (\ref{8ii}), we conclude
$$
\|F_{2N}^\textsc{e}(q)-\be_k({\tt P}_k(q))\|<\sqrt{\left(2s+\frac{\cte}{\sqrt{N}} \right)^2+\left(r+\frac{\cte}{N} \right)^2}+\frac{\cte}{N}\;.
$$
So, for a large enough $N$, it is obvious that $\|F_{2N}^\textsc{e}(q)-\be_k({\tt P}_k(q))\|<R$ in this last case. 

This completes the proof of (d3.4) and concludes the proposition.
\end{proof}

From the above proposition  it is straightforward to check that $\widetilde{X}=F_{2N}^\textsc{e}:\overline{M(\widetilde{\mathcal{ J}})}\to \R^3$  proves Lemma \ref{nadira}.




\subsection{Proof of Lemma \ref{main}}

Consider $c_0,$ $r_1$ and $\rho_1$ three positive constants to be specified later, and define
$$
r_n=\sqrt{r_{n-1}^2+\left(\frac{2c_0}{n}\right)^2}+\frac{c_0}{n^2}\quad\text{ and }\quad \rho_n=\rho_1+\sum_{i=2}^n \frac{c_0}{i}\,,\quad \forall\,n\geq 2\;.
$$
The constants $r_1$ and $c_0$ have to be chosen so that
\begin{equation}\label{lim-r}
\lim_{n\to\infty} r_n<\frac{\mu}{2}\;.
\end{equation}

In order to apply Lemma \ref{nadira}, we consider a family of analytic cycles in $M'$, $\Sigma_i$, $i=1, \ldots,\textsc{e}$, such that $\g_{i}\subset \mathcal{T}(\Sigma_i)$, for $i=1, \ldots,\textsc{e}$, where $\mathcal{T}(\Sigma_i)$ is a tubular neighborhood of the curve $\Sigma_i$ described at the beginning of this section.
 
Hereafter, we will construct a sequence $\chi_n=\{\mathcal{ J}_n,X_n,\ep_n\}$ consisting of:
\begin{itemize}
\item $\mathcal{ J}_n=\{\g_{n,1},\ldots,\g_{n,\textsc{e}}\}$ is a multicycle with $\g_{n,i}\subset \mathcal{T}(\Sigma_i)$ for $i=1,\ldots, \textsc{e}$;

\item $X_n:\overline{M(\mathcal{ J}_n)}\to \R^3$ is a conformal minimal immersion;

\item $\{\ep_n\}$ is a decreasing sequence of positive real numbers with $\ep_n<c_0/n^2.$
\end{itemize}

\begin{claim}\label{Af-Meeks}
The sequence $\{\chi_n\}$ can be constructed to satisfy:
\begin{enumerate}[\rm (A$_n$)]
\item $\mathcal{J}^\mu<\mathcal{ J}_{n-1}^{\ep_n}<\mathcal{ J}_n<\mathcal{ J}_{n-1};$

\item $\dist_{(\overline{M(\mathcal{ J}_n)},X_n)}(p,\mathcal{ J}_n)>\rho_n$, for all $p \in \mathcal{J}^\mu;$

\item $\|X_n(p)-X({\tt P}_k(p))\|<r_n,$ $\forall p\in \g_{n,k},$ $\forall k=1,\ldots,\textsc{e}.$
\end{enumerate}
\end{claim}
\noindent Notice that (A$_n$) only holds for $n\geq2.$ Once again, the sequence will be obtained following a inductive method. For the first term, we choose $X_1=X$ and $\mathcal{ J}_1=\mathcal{ J}.$ Finally, we take $\rho_1$ and $\ep_1$ satisfying
$$
\rho_1<\dist_{(X_1,\overline{M(\mathcal{ J}_1)})}(p,\mathcal{ J}_1), \; \, \mbox{for all $p \in \mathcal{J}^\mu$}\quad \text{ and }\quad \ep_1<\min\{c_0,r_1\}\;.
$$
Moreover, we take $\ep_1$ small enough so that $\gamma_i^{\ep_1} \subset \mathcal{T}(\Sigma_i)$, $i=1, \ldots, \textsc{e}$, and 
\begin{equation}\label{entornito}
\|X(p)-X({\tt P}_k(p))\|<r_1  <\frac{\mu}{2} \;,
\end{equation}
for any $p$ in the connected component of $\overline{M(\mathcal{ J})}-M(\mathcal{ J}^{\ep_1})$ around $\g_k,$ $\forall k=1,\ldots,\textsc{e}.$

Assume now that we have constructed $\chi_1,\ldots, \chi_{n-1}.$ In order to define $\chi_n$ we take a real number $\ep_n<\min \{\ep_{n-1},\frac{c_0}{n^2}\}.$ Then we consider the multicycle $\mathcal{ J}_n$ and the immersion $X_n:\overline{M(\mathcal{ J}_n)}\to\R^3$ given by Lemma \ref{nadira}, for the data
$$
X=X_{n-1}\;,\quad \mathcal{ J}=\mathcal{ J}_{n-1}\;,\quad r=r_{n-1}\;,\quad s=\frac{c_0}{n}\quad \text{ and }\quad \ep=\ep_n\;.
$$
So, we get $\chi_n$ satisfying properties (A$_n$), (B$_n$) and (C$_n$).

From (A$_n$), (B$_n$) and the fact that the sequence $\{\rho_n\}_{n\in\N}$ diverges, we find $n_0\in \N$ with $\dist_{(\overline{M(\mathcal{ J}_n)},X_n)}(p,\mathcal{ J}_n)>\l,$ $\forall p\in \mathcal{ J}^\mu,$ $\forall n\geq n_0.$ Choose $\widehat{X}=X_{n_0}$ and $\widehat{\mathcal{ J}}=\mathcal{ J}_{n_0}.$ Properties (L4.a) and (L4.b) trivially hold. Now, taking (\ref{entornito}), (C$_{n_0}$) and  (\ref{lim-r}) into account, we obtain
$$
\|X(p)-\widehat{X}(p)\|\leq \|X(p)-X({\tt P}_{k}(p))\|+\|X({\tt P}_{k}(p))-\widehat{X}(p)\| <
\frac{\mu}{2}+r_{n_0}<\mu\;,\quad \forall p\in \widehat{\g}_k\;,\; \; \forall k=1,\ldots,\textsc{e}\;.
$$
Hence, $\|X(p)-\widehat{X}(p)\|\leq \mu$ for any $p\in \widehat{\mathcal{ J}}.$ Finally, the Maximum Principle guarantees that this inequality occurs for any $p\in M(\widehat{\mathcal{ J}}),$ so we have checked (L4.c).

\begin{remark} \label{re:delta}
From the arguments of the above proof, it is almost trivial to deduce that:
\[ \delta^H\left(X ( \overline{M(\mJ)}), \widehat X (\overline{ M(\widehat \mJ)})\right)<2 \, \mu.\]
This estimation will be important to prove Theorem \ref{teo-ultimo}.
\end{remark}




\section{Joining together properness and completeness}
As the title indicates, in this section we put together the information obtained in the previous two sections in order to state the precise lemma that we will use in the proof of the main theorems.

\begin{lemma}\label{3+5} Let $\mJ$ be a multicycle in $M$, $p_0 \in M(\mJ)$,  and 
$X:\overline{M(\mathcal{ J})}\to \R^3$  a conformal minimal immersion with $X(p_0)=0$. Consider $E$ and $E'$ bounded convex regular domains, with $0\in E\subset \overline{E}\subset E'$,   and  let $a$ and $\ep$ be positive constants satisfying that 
$p_0 \in M(\mathcal{J}^\ep)$ and 
\begin{equation}\label{al}
X(\overline{M(\mathcal{ J})}-M(\mathcal{ J}^\ep))\subset E-\overline{E_{-a}}\;.
\end{equation}
Then, for any $b>0$ there exist a multicycle $\widetilde{\mathcal{ J}}$ and a conformal minimal immersion $Y:\overline{M(\widetilde{\mathcal{ J}})}\to\R^3$ such that $Y(p_0)=0$ and
\begin{enumerate}[\rm (L\ref{3+5}.a)]
\item $\mathcal{ J}^\ep<\widetilde{\mathcal{ J}}<\mathcal{ J};$

\item $\dist_{(\overline{M(\widetilde{\mathcal{ J}})},Y)}(p,\mathcal{ J}^\ep)>1/\ep,$ $\forall p\in \widetilde{\mathcal{J}};$

\item $Y(\widetilde{\mathcal{ J}})\subset E'-\overline{E'_{-b}};$

\item $Y(\overline{M(\widetilde{\mathcal{ J}})}-M(\mathcal{ J}^\ep))\subset \R^3-E_{-2b-a};$

\item $\|X-Y\|<\ep$ in $M(\mathcal{ J}^\ep).$

\end{enumerate}
Furthermore if $E$ is {strictly convex}, the immersion $Y$ also satisfies:
\begin{enumerate}[{\rm (L5.f)}]
\item $\|X-Y\|<\mathsf{m}(a,b,\ep,E,E'):= \ep+\sqrt{ \frac{2 (\delta^H(E,E')+a+2 b)}{\kappa_1(\partial E)}+(\delta^H(E,E')+a)^2},$ in $M(\widetilde{\mathcal{ J}}).$ 
\end{enumerate}
\end{lemma}


\begin{proof}
First, we apply Lemma \ref{main} to the immersion $X$, for $\l>1/\ep$ and a small enough  $\mu>0$  which will be determined later. Then, we get a new multicycle $\widehat{\mathcal{J}}$ and a immersion $\widehat X: \overline{M(\widehat{ \mathcal{J}})} \rightarrow \R^3$, such that:
\begin{enumerate}[(a)]
\item $\mathcal{J}^\ep < \mJ^\mu<\widehat \mJ<\mJ;$
\item $\dist_{(\overline{M(\widehat \mJ)}, \widehat X)}(\widehat \mJ, \mJ^\mu)> \lambda;$
\item $ \|X-\widehat X\|<\mu,$ in $M(\widehat \mJ).$
\end{enumerate}
If $\mu$ is sufficiently small, then $\widehat X(\widehat \mJ) \subset E-\overline{E_{-a}}.$ Thus, we can find $\nu>0$ so that $ \mJ^\mu <\widehat \mJ^{\nu} $ and:
\begin{gather} \label{tack1} \dist_{(\overline{M(\widehat \mJ^\nu)}, \widehat X)}(\widehat \mJ^\nu, \mJ^\mu)> \lambda \\
\widehat X \left(\overline{M(\widehat \mJ)}-M(\widehat \mJ^{\nu})\right)\subset E-\overline{E_{-a}}.
\end{gather}
At this point, we apply Lemma \ref{properness} to the following data:
$$X=\widehat X, \quad E=E_{-a}, \quad E', \quad \mJ_0=\widehat \mJ, \quad \mJ'=\widehat \mJ^{\nu}, \quad b_2=b, $$
and arbitrary $b_1>0.$ Hence, we obtain a new multicycle $\widetilde \mJ$, $\widehat \mJ^{\nu}<\widetilde \mJ<\widehat \mJ$, and a minimal immersion $Y:\overline{M(\widetilde \mJ)} \rightarrow \R^3, $ satisfying:
\begin{enumerate}[(A)]
\item $\|Y- \widehat X\|<b_1,$ in $M(\widehat \mJ^\nu);$
\item $Y(\widetilde \mJ) \subset E' -\overline{E'_{-b}}; $
\item $Y( \overline{M(\widetilde \mJ)}-M(\widehat \mJ^\nu)) \subset \R^3-E_{-2 b-a}. $
\end{enumerate}
Furthermore, if $E$ is strictly convex, then we have the extra information provided by Remark \ref{k1>0}:
\begin{enumerate}[{\rm(D)}]
\item $\|Y-\widehat X\|<\mathcal{M}(b,E_{-a},E')<\sqrt{ \frac{2 (\delta^H(E,E')+a+2 \,b)}{\kappa_1(\partial E)}+(\delta^H(E,E')+a)^2},\;$ in $\overline{M(\widetilde \mJ)}-M(\widehat \mJ^\nu).$
\end{enumerate}
Item (A) says to us that $Y$ converges to $\widehat X$ uniformly on $M(\widehat \mJ^\nu)$, as $b_1 \to 0.$ Therefore, if $b_1$ is small enough we also have $\dist_{(M(\widehat \mJ^\nu), Y)}(\widehat \mJ^\nu, \mJ^\mu)> \lambda$ (see \eqref{tack1}), which implies (L5.b). Item (L5.c) directly follows from (B). Moreover, (c) and (A) give (L5.e), provided that $\mu+b_1<\ep.$

Taking \eqref{al}, (c) and (A) into account, we can deduce that $Y( \overline{M(\widehat \mJ^\nu)}-M(\mJ^\ep)) \subset \R^3-\overline{E_{-a}}$, provided that $\mu$ and $b_1$ are sufficiently small. So, the above inclusion and (C) demonstrate (L5.d).

Finally, if $E$ is strictly convex, then (c), (A) and  (D) imply (L5.f), provided that $\mu$ and $b_1$ are small enough.
\end{proof}




\section{Density theorems for complete minimal surfaces in $\R^3$} \label{sec:teoremas}

Now, we are able to prove the theorems stated in the introduction. Although all the theorems of this section are stated in terms of Riemann surfaces with boundary that are open regions of compact Riemann surfaces, this does not represent any restrictions over our work. In order to prove Theorem A in the introduction we notice that any Riemann surface with finite topology and analytic boundary can be seen as the closure of an open region of a compact Riemann surface (see \cite{a-s}.) 

\begin{remark} \label{alfas}
In this section, we will use several times the sequence of positive reals given by:
$$ \a_1:=\tfrac 12 \;{\rm e}^{1/2}, \quad \a_n:= {\rm e}^{-1/2^n}, \; \, \mbox{\rm for $n>1$.} $$
Notice that $0<\alpha_i<1$ and $\{ \prod ^n_{i=1} \alpha_i\}_{n \in \N}$ converges to $1/2.$
\end{remark}
\begin{theorem} \label{teo-A}
Let $D$ and $D'$ be two bounded, convex regular domains satisfying $0 \in \overline{D} \subset D'.$ Let $\varphi: \overline{M(\Gamma)} \rightarrow \R^3$ be a conformal minimal immersion, where $\Gamma$ is a multicycle in $M'$. Assume that $\varphi(p_0)=0$ and $\varphi(\Gamma) \subset D-\overline{D_{-d}}$ where $p_0 $ is a point in $ M(\Gamma)$ and $d$ is a positive constant. 

Then for any $\mu>0$, there exists a domain $M_\mu$ in $M'$, with $\overline{M(\Gamma)} \subset M_\mu$ and there exists a complete proper minimal immersion $\varphi_\mu: M_\mu \rightarrow D'$ such that:
\begin{enumerate}[\rm (a)]
\item $\| \varphi_\mu-\varphi \| <\mu$ in $M(\Gamma)$;
\item $\varphi_\mu(M_\mu - M(\Gamma)) \subset D' - \overline{D_{-2 d-\mu}}.$
\end{enumerate}
\end{theorem}
\begin{proof}
First of all, we define a sequence $\{E^n \}$ of bounded convex regular domains in the following way. Consider $\nu>0$ small enough to satisfy that $D'_{-\widetilde\nu}$ exists, $\overline{D}\subset D'_{-\widetilde\nu}$, 
where $\widetilde\nu=\sum_{k=2}^\infty \nu/k^2$. Then, we define 
$$E^1:= D \quad \mbox{and} \quad E^n:= D'_{-\sum_{k=n}^\infty \nu/k^2}, \; \; n\geq 2.$$ 

We also take a decreasing sequence of positive reals $\{b_n\}$ with  $b_1=d$, and: 
$$b_n<\min\left\{\frac{\dist_{\R^3}(\partial E^n, \partial E^{n+1})}2, \frac{d}{2}\right\}, \quad \mbox{for $n>1$}.$$

Next, we use Lemma \ref{3+5} to construct a sequence 
$$\chi_n=(\varphi_n:\overline{M(\Gamma_n}) \rightarrow \R^3,\Gamma_n, \varepsilon_n,\xi_n),$$
where $\varphi_n$ are conformal minimal immersions with $\varphi_n(p_0)=0$, $\Gamma_n$ are multicycles,  and $ \{     \varepsilon_n\}$, $\{\xi_n\}$  are sequences of positive numbers decreasing to zero, and satisfying $\sum_{k=1}^\infty \varepsilon_k <\mu.$

Furthermore, the sequence $\varphi_n:\overline{M( \Gamma_n)} \rightarrow \R^3$ must satisfy the following properties:
\begin{enumerate}[(A$_{n}$)]
\item $\Gamma_{n-1}^{\xi_{n-1}}< \Gamma_{n-1}^{\varepsilon_n}< \Gamma_n^{\xi_n}< \Gamma_n< \Gamma_{n-1}$;
\item $\| \varphi_n(p)-\varphi_{n-1}(p)\| < \varepsilon_n$, $\forall p \in M(\Gamma_{n-1}^{\varepsilon_{n}})$;
\item $\metri{\varphi_n}(p) \geq \alpha_n \cdot \metri{\varphi_{n-1}}(p)$, $\forall p \in M( \Gamma_{n-1}^{\xi_{n-1}})$, where $\{\alpha_i\}_{i \in \N}$ is given by Remark \ref{alfas} (recall that $\metri{\varphi_n}$ means the Riemannian metric induced by $\vp_n$);
\item $ 1/\varepsilon_n < \dist_{(\overline{M( \Gamma_n^{\xi_n})},\varphi_n)}(\Gamma_{n-1}^{\xi_{n-1}},\Gamma_n^{\xi_n})$;
\item $\varphi_n(p)\in E^{n}- \overline{(E^{n})_{-b_n}}$, for all $p \in \Gamma_n$;
\item $\varphi_n(p)\in \R^3- (E^{n-1})_{-b_{n-1}-2b_n}$, for all $p \in \overline{M( \Gamma_n)}- M( \Gamma_{n-1}^{\varepsilon_n})$.
\end{enumerate}

The sequence $\{\chi_n\}$ is constructed in a recursive way. To define $\chi_1$, we take $\varphi_1:=\varphi$ and $\xi_1>0$ small enough so that $\Gamma^{-\xi_1}$ is well-defined, $\varphi$ is defined in $\overline{M(\Gamma^{-\xi_1})}$  and 
\begin{equation}\label{eq:dirigir}
\varphi \left(\overline{M(\Gamma^{-\xi_1})}-M(\Gamma)\right) \subset D-\overline{D_{-d}}. 
\end{equation}
By definition $\Gamma_1:=\Gamma^{-\xi_1}.$ In particular Property (E$_1$) holds. The other properties do not make sense for $n=1$.

Suppose that we have $\chi_1, \ldots, \chi_n$. In order to construct $\chi_{n+1}$, we consider the following data:
$$E=E^n,\quad E'=E^{n+1},\quad a=b_{n}, \quad X=\varphi_n,\quad \mathcal{J}=\Gamma_n.$$
Furthermore, Property (E$_n$) tells us that $X(\mathcal{J})\subset E - E_{-a}$. Then it is straightforward that we can find a small enough positive constant $\varkappa$, such that Lemma \ref{3+5} can be applied to the aforementioned data, and for any $\epsilon \in ]0,\varkappa[$. 

Take a sequence $\{\widehat{\varepsilon}_m \}_{m\in \N}\searrow 0$, with $\widehat{\varepsilon}_m<\min\{\varkappa , b_{n+1}\}$, $\forall m$. For each $m$, we consider $\mathcal{J}'_m$ and $Y_m: \overline{M( \mathcal{J}'_m)} \rightarrow \R^3$  given by Lemma \ref{3+5}, for the above data and $\epsilon=b=\widehat{\varepsilon}_m$.
If $m$ is large enough, Assertions (L5.a) and (L5.e) in Lemma \ref{3+5} tell us that $\Gamma_n^{\xi_n} <\mathcal{J}'_m$ and the sequence $\{Y_m\}$ converges to $\varphi_n$ uniformly in $\overline{ M( \Gamma_n^{\xi_{n}})}$. In particular, $\{\metri{Y_m}\}_{m \in \N}$ converges uniformly to $\metri{\varphi_n}$  in $\overline{ M( \Gamma_n^{\xi_{n}})}$. Therefore there is a $m_0 \in \N$ such that:
\begin{eqnarray} 
\Gamma_n^{\xi_n} < &  \Gamma_n^{\hat \varepsilon_{m_0}} & <  \mathcal{J}'_{m_0}, \label{clavao1}\\
\metri{Y_{m_0}}& \geq & \alpha_{n+1} \cdot \metri{\varphi_n} \qquad \hbox{in } M( \Gamma_n^{\xi_{n}}).\label{lambdas}
\end{eqnarray}
We define $\varphi_{n+1}:=Y_{m_0}$, $\Gamma_{n+1}:=\mathcal{J}'_{m_0}$, and $\varepsilon_{n+1}:=\widehat{\varepsilon}_{m_0}$.  From (\ref{clavao1}) and Statement (L5.b), we infer that $1/\varepsilon_{n+1} < \dist_{\left(\overline{ M( \Gamma_{n+1})}, \varphi_{n+1} \right)}(\Gamma_n^{\xi_n},\Gamma_{n+1})$. Finally, take $\xi_{n+1}$ small enough such that (A$_{n+1}$) and (D$_{n+1}$) hold. 
The remaining properties directly follow from (\ref{clavao1}), (\ref{lambdas}) and Lemma \ref{3+5}. This concludes the construction of the  sequence $\{\chi_n\}_{n \in \N}$.
\bigskip

Now, we extract some information from the properties of $\{\chi_n\}$. First, from (B$_n$), we deduce that $\{\varphi_n\}$ is a Cauchy sequence, uniformly on compact sets of $M_\mu=\bigcup_n M( \Gamma_n^{\varepsilon_{n+1}})=\bigcup_n M( \Gamma_n^{\xi_{n}})$, and so $\{\varphi_n\}$ converges on $M_\mu$. If one employs the properties (A$_n$), then the set $M_\mu$ is an expansive union of domains with the same topological type as $M(\Gamma)$. Therefore, elementary topological arguments give us that $M_\mu$ has the same topological type as $M(\Gamma)$. Let $\varphi_\mu:M_\mu\rightarrow \R^3$ be the limit of $\{\varphi_n\}$. Then $\varphi_\mu$ has the following properties:
\begin{itemize}
\item $\varphi_\mu$ is a conformal minimal immersion, (Properties  (C$_n$));
\item $\varphi_\mu:M_\mu \longrightarrow D'$ is proper. Indeed, consider a compact subset $K \subset D'$. Let $n_0$ be a natural so that $$K \subset (E^{n-1})_{-b_{n-1}-2 b_n-\sum_{k\geq n} \varepsilon_k}, \quad \forall n  \geq n_0.$$ From Properties (F$_n$), we have $\varphi_n(p)\in \R^3- (E^{n-1})_{-b_{n-1}-2 b_n}$, $\forall p \in \overline{M( \Gamma_n)}- M( \Gamma_{n-1}^{\varepsilon_n})$. Moreover, taking into account (B$_k$), for $k \geq n$, we obtain 
\begin{equation} \label{violeta}
\varphi_\mu(\overline{M( \Gamma_n)}- M( \Gamma_{n-1}^{\varepsilon_n}))\subset \R^3- (E^{n-1})_{-b_{n-1}-2 b_n-\sum_{k\geq n} \varepsilon_k}.
\end{equation}
Then, we have $\varphi_\mu^{-1}(K)\cap \left(\overline{M( \Gamma_n)} - M( \Gamma_{n-1}^{\varepsilon_n})\right)=\emptyset$ for $n\geq n_0$. This implies that $\varphi_\mu^{-1}(K)\subset M( \Gamma_{n_0-1}^{\varepsilon_{n_0}})$, and so it is compact in $M_\mu$.
\item Completeness of $\varphi_\mu$ follows from Properties (D$_n$), (C$_n$), and the fact that $\{1/\varepsilon_n\}_{n \in \N}$ diverges.
\item Statement {(a)} in the theorem is a direct consequence of Properties (B$_n$) and the fact $\sum_{n=1}^\infty \varepsilon_n <\mu$.
\item In order to prove Statement {(b)}, we consider $p \in M_\mu- M( \Gamma)$. If there exists  $n \in \N$ such that $p \in \overline{M( \Gamma_n)} - M(\Gamma_{n-1}^{\varepsilon_n}),$ then \eqref{violeta} implies
$\varphi_\mu(p) \in \R^3 - D_{-2 d-\mu}.$
If $p \in M(\Gamma_1^{\varepsilon_2})-M(\Gamma_1)$, then we use properties (B$_k$), $k \geq 1$, and \eqref{eq:dirigir}
to obtain $$\varphi_\mu(p) \in \R^3 - D_{-d-\sum_{k\geq 1} \varepsilon_k}\subset \R^3 - D_{-2 d-\mu}.$$
\end{itemize}
\end{proof}


If we follow the proof of the above theorem, but making use of Lemma \ref{main} and Remark \ref{re:delta} instead of Lemma \ref{3+5}, then we obtain the following theorem:
\begin{theorem}\label{teo-ultimo}
Let $\varphi: \overline{M(\Gamma)} \rightarrow \R^3$ be a conformal minimal immersion, where $\Gamma$ is a multicycle in $M'$.  Then for any $\mu>0$, there exists a domain $M_\mu$ in $M'$, with $\overline{M(\Gamma)} \subset M_\mu$ and there exists a complete proper minimal immersion $\varphi_\mu: M_\mu \rightarrow \R^3$ such that:
\begin{enumerate}[\rm (a)]
\item $\| \varphi_\mu-\varphi \| < \mu$ in $M(\Gamma)$;
\item $\de^H\left(\vp(\overline{M(\Gamma)}),\overline{\vp_\mu(M_\mu)}\right)<2 \, \mu.$
\end{enumerate}
\end{theorem}

Under the assumption of strictly convexity we can sharpen the previous arguments in order to prove the following theorem.

\begin{theorem}\label{teo-B}
Let $C$ be a strictly convex bounded regular domain of $\R^3$. Consider a multicycle $\Gamma$ in the Riemann surface $M'$ and $\vp:\overline{M(\Gamma)}\to \overline{C}$ a conformal minimal immersion satisfying $\vp(\Gamma)\subset \partial C$. 

Then, for any $\ep>0$, there exist a subdomain, $M_\ep,$ with the same topological type as $M(\Gamma)$, $\overline{M(\Gamma^\ep)}\subset M_\ep \subset \overline{M_\ep}\subset M(\Gamma)$, and a complete proper conformal minimal immersion $\vp_\ep:M_\ep \to C$ so that 
\begin{equation*} \label{rafa}
\|\vp-\vp_\ep\|<\ep\;,\quad \text{in }M_\ep\;.
\end{equation*}
\end{theorem}


\begin{proof}

Consider $t_0>0$ so that, for any $t\in]-t_0,0[$, we have:
\begin{itemize}
\item  $C_t$ is well-defined;
\item $\Gamma_t:=\vp^{-1}(\partial C_t \cap \vp(M(\Gamma)))$ is a multicycle satisfying $\Gamma^\ep<\Gamma_t . $
\end{itemize}

Fix $c_1>0$ small enough so that 
$\sum_{k\geq 1}{c_1^2}/{k^4}<\min\{t_0,\ep\}.$ 
At this point, for any $n\geq 1,$ consider a positive constant $t_n=\sum_{k\geq n}c_1^2/k^4$ and a strictly convex bounded regular domain $E^n=C_{-t_n}.$ We also take a decreasing sequence of positive reals $\{b_n\}_{n\in\N},$ with $b_n<c_1^2/n^4,$ $\forall n\in\N.$

Now, we use Lemma \ref{3+5} to construct, for any $n\in\N$, a family $\chi_n=\{\mathcal{ J}_n,X_n,\ep_n,\xi_n\}\;,$ where
\begin{itemize}
\item $\mathcal{ J}_n$ is a multicycle;

\item $X_n:\overline{M(\mathcal{ J}_n)}\to C$ is a conformal minimal immersion;

\item $\{\ep_n\}_{n\in\N}$  and $\{\xi_n\}_{n\in\N}$ are  sequences of positive real numbers converging to zero and satisfying $\ep_n<c_1/n^2$.
\end{itemize}

Notice that the function given in (L5.f) satisfies $
\mathsf{m}(b_n,b_{n+1},\ep_{n+1},E^n,E^{n+1})<\frac{c_1}{n^2}\left( 1+2 \sqrt{\frac{{c_1}^2}{n^4}+\frac{2}{\kappa_1(\partial C)}}\, \right),$ $\forall n\in\N,$
therefore, we can choose $c_1$ sufficiently small so that
\begin{equation}\label{cota-M}
\sum_{n=2}^\infty \mathsf{m}(b_{n-1},b_n,\ep_n,E^{n-1},E^n)<\ep\;.
\end{equation}
We will construct the sequence $\{\chi_n\}_{n\in\N}$ so that the following properties hold:
\begin{enumerate}[(A$_{n}$)]
\item $\mJ_{n-1}^{\xi_{n-1}}< \mJ_{n-1}^{\ep_n}< \mJ_n^{\xi_n}< \mJ_n< \mJ_{n-1}$;

\item $\| X_n(p)-X_{n-1}(p)\| < \ep_n$, $\forall p \in M( \mJ_{n-1}^{\ep_{n}})$;
\item $\metri{X_n}(p) \geq \alpha_n \cdot \metri{X_{n-1}}(p)$, $\forall p \in M( \mJ_{n-1}^{\xi_{n-1}})$, where $\{\alpha_i\}_{i \in \N}$ is given by Remark \ref{alfas};
\item $ 1/\ep_n < \dist_{(\overline{M( \mJ_n^{\xi_n})},X_n)}(\mJ_{n-1}^{\xi_{n-1}},\mJ_n^{\xi_n})$;

\item $X_n(p)\in E^{n}- \overline{(E^{n})_{-b_n}}$, for all $p \in \mJ_n$;
\item $X_n(p)\in \R^3- (E^{n-1})_{-b_{n-1}-2 b_n}$, for all $p \in \overline{M( \mJ_n)}- M( \mJ_{n-1}^{\ep_n})$;
\item $\|X_n-X_{n-1}\|<\mathsf{m}(b_{n-1},b_n,\ep_n,E^{n-1},E^n)$ in $M(\mathcal{ J}_n).$
\end{enumerate}

The construction of the sequence $\{\chi_n\}_{n\in \N},$ is as in the proof of Theorem \ref{teo-A}, except for properties (G$_n$) that are consequence of the successive use of (L5.f) in Lemma \ref{3+5}. To define $\chi_1,$ we take $X_1=\vp,$ $\mathcal{ J}_1=\Gamma_{t_1}$ and  appropriate $\ep_1$ and $b_1.$ We choose $\xi_1$ so that $\Gamma^\ep<\mJ_1^{\xi_1}.$
Observe that in this case properties (G$_n$), $n \in \N$, and (\ref{cota-M}) guarantee that $\|\vp_\ep -\vp \| < \ep,$ in $M_\ep$.
\end{proof}

\begin{corollary}\label{gordo}
Let $D'$  be  a convex domain (not necessarily bounded or smooth) in $\R^3$. Consider $\mathcal{J}$ a multicycle in $M'$ and $\vp:\overline{M(\mathcal{J})} \longrightarrow \R^3$ a conformal minimal immersion satisfying:
\begin{equation}
\vp(\mathcal{J}) \subset D - \overline{D_{-d}}.
\end{equation}
where $D$ is a bounded convex regular domain satisfying $\overline{D} \subset D'$ and $d>0$ is a constant.

Then, for any $\varepsilon>0$, there exist a subdomain, $M_\varepsilon,$ with $\overline{M(\mathcal{J})}\subset M_\varepsilon \subset \overline{M_\varepsilon}\subset M(\mathcal{J}^{-\varepsilon}),$  and a complete proper conformal minimal immersion $\vp_\varepsilon:M_\varepsilon \longrightarrow D'$ so that
\begin{enumerate}[\rm (A)]
\item $\|\vp-\vp_\varepsilon\|<\varepsilon\;,\quad \text{in } M(\mathcal{J}) \;;$
\item $\vp_\varepsilon(M_\varepsilon -M(\mathcal{J})) \subset \R^3-D_{-2 d-\varepsilon}\;.$
\end{enumerate}
\end{corollary}
\begin{proof}
Without loss of generality we can assume $0\in D$ and $\varphi(p_0)=0$, for a certain $p_0$ in $M(\mathcal{J}).$ Define $V^n:= \left[(1-1/n) \cdot D' \right]\cap \B(0,n),$ where by $(1-1/n) \cdot D'$ we mean the set $\{(1-1/n)\cdot x \; | \; x \in D'\}.$ Making use of Minkowski's theorem (see page \pageref{th:minko}) we can guarantee, for each $n \in \N$,  the existence of a regular bounded domain $\widehat V^n$ in $\R^3$ such that $\overline{V^n} \subset \widehat V^n\subset \overline{\widehat V^n} \subset V^{n+1}.$ Notice that $\{\widehat V^n\}_{n \in \N}$ is an expansive sequence of bounded convex regular  domains whose limit is $D'.$ Then there exists $k \in \N$ so that $\overline{D} \subset \widehat V^m$ for any $m \geq k.$ Taking all these arguments into account, we define the following sequence of open convex domains:
$$E^1 := D, \quad E^n:= \widehat V^{n+k-2}, \; n \geq 2.$$

\noindent Following the scheme of the previous proofs, we will construct a sequence $\Xi_n:=(\mathcal{J}_n, \varphi_n:\overline{M(\mJ_n)} \to \R^3, d_n, \varepsilon_n),$
where:
\begin{itemize}
\item $\mJ_n$ is a multicycle in $M'$;
\item $\varphi_n:\overline{M(\mJ_n)} \to \R^3$ is a conformal minimal immersion;
\item $\{d_n\}_{n \in \N}$ and $\{\varepsilon_n\}_{n \in \N}$ are two sequences of positive real numbers  decreasing to $0$. Moreover we want that $ \sum_{k =1}^\infty \varepsilon_k< \varepsilon, $ to do this we will choose $ \varepsilon_n <\frac{6 \varepsilon}{\pi^2 n^2},$ $n \in \N.$
\end{itemize}
The limit of  $\{\Xi_n\}_{n \in \N}$ will provide the minimal immersion we are looking for. To do this we need that $\Xi_n$ satisfies the following properties:

\begin{enumerate}[(I$_n$)]
\item $\mJ_{n-1}<\mJ_n$;
\item $\| \vp_n(p)-\vp_{n-1}(p)\| < \varepsilon_n$, for all $p$ in $M(\mJ_{n-1})$;
\item $\metri{\varphi_n} \geq \alpha_n \cdot \, \metri{\vp_{n-1}}$ in $\overline{M(\mJ_{n-1})}$, where 
 $\{\a_i\}_{i\in\N}$ is the sequence of Remark \ref{alfas};
\item $\dist_{(\overline{M(\mJ_n)},\vp_n)}(p_0, \mJ_n)>n-1;$
\item $\vp_n(\mJ_n) \subset E^n-\overline{(E^n)_{-d_n}};$
\item $\vp_n \left( \overline{M(\mJ_n)} - M(\mJ_{n-1}) \right) \subset \R^3-E^{n-1}_{-2 d_{n-1}-\varepsilon_n}.$
\end{enumerate}

Once again the  sequence $\{\Xi_n\}_{n \in \N}$ satisfying the above properties is defined following an inductive process. The elements of $\Xi_1$ are $\vp_1:=\vp$, $d_1=d$, $\mJ_1=\mJ$ and $\varepsilon_1<\frac{6 \, \varepsilon}{\pi^2}.$ 

Assume now we have defined $\Xi_n$. To construct the element $\Xi_{n+1}$ we apply Theorem \ref{teo-A} to the minimal immersion $\vp_n:M(\mJ_n) \to E^n \subset E^{n+1},$ where $E^n$, $E^{n+1}$, $\varepsilon_{n+1}$, and $d_n$ play the role of $D$, $D'$, $\mu$, and $d$ in the statement of Theorem \ref{teo-A}, respectively. Then we get a domain $M_{\varepsilon_{n+1}}$ in $M'$, with $\overline{M(\Gamma_n)}\subset M_{\varepsilon_{n+1}}$, and a complete proper minimal immersion $\vp_{n+1}: M_{\varepsilon_{n+1}} \to E^{n+1}$ satisfying: 
\begin{eqnarray}
\|\vp_{n+1}-\vp_n \| & < & \varepsilon_{n+1}, \; \mbox{in } \; M(\mJ_n); \label{mari} \\
\vp_{n+1} (M_{\varepsilon_{n+1}}-M(\mJ_n)) & \subset & \R^3-E^{n}_{-2 d_{n}-\varepsilon_{n+1}}.\label{quita}
\end{eqnarray}
From \eqref{mari} we have that (III$_{n+1}$) holds provided that $\varepsilon_{n+1}$ is taken small enough. As $\vp_{n+1}$ is complete and proper, then it is possible to find $\mJ_{n+1}$ satisfying (I$_{n+1}$), (IV$_{n+1}$) and (V$_{n+1}$).  Properties (II$_{n+1}$) and (VI$_{n+1}$) are consequence of \eqref{mari} and \eqref{quita}, respectively.

At this point we define $M_{\varepsilon}:= \cup_{n=1}^\infty M(\mJ_n)$ and $\vp_\varepsilon: M_\varepsilon \to D'$ as the uniform limit of the sequence $\{ \vp_n \}_{n \in \N}.$ Following similar arguments to those used in the proof of Theorem \ref{teo-A}, it is easy to check that $\vp_\varepsilon$ is the minimal immersion that proves the corollary.
\end{proof}

\section{The construction of a complete proper minimal surface with uncountably many ends} 

\label{sec:uncountable}
The most interesting application of the results in the preceding  section is the construction of  the first examples of complete properly immersed minimal surfaces in Euclidean space with an uncountable number of ends. It is important to note  that this kind of surfaces cannot be embedded as a consequence of a result by Collin, Kusner, Meeks and Rosenberg \cite{CKMR}. Given $p\in \C$ and $r>0$, we will write $\D(p,r)=\{z \in \C \; | \; |z-p|<r\}.$ As usual,  the unit disk will be denoted by $\D$, instead of $\D(0,1).$
\begin{figure}
	\begin{center}
		\includegraphics[width=0.55\textwidth]{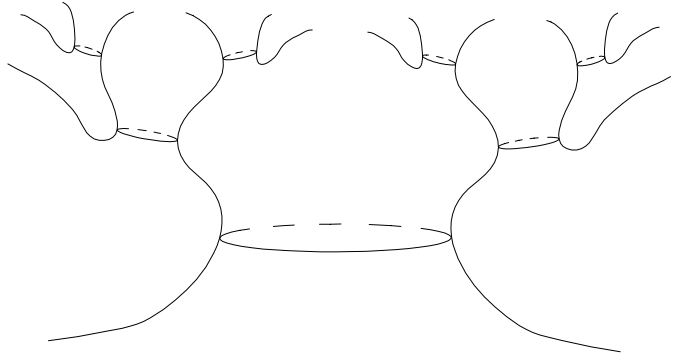}
	\end{center}
	\caption{The construction of a minimal surface with uncountably many ends consists of modifying a given minimal surface in $\R^3$ by adding more and more ends. In each step of this procedure we add two new ends in a neighborhood of each end of the previous surface.}
	\label{fig:catenoides1}
\end{figure}
\begin{figure}[ht]
\begin{center}
\scalebox{0.35}{\includegraphics{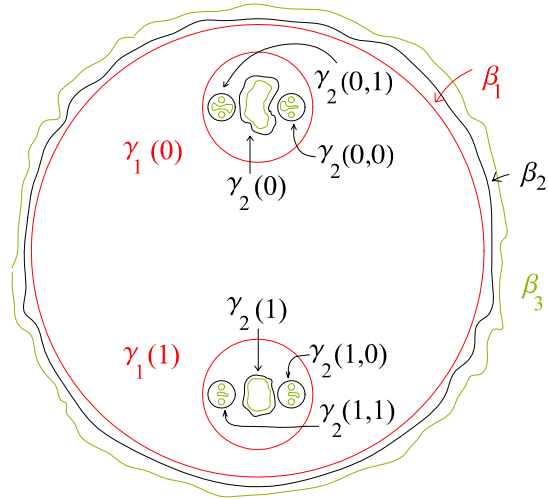}} 
\end{center}
\caption{The multicycles $\Gamma_1$, $\Gamma_2$, and $\Gamma_3$.}\label{quesito}
\end{figure}
\begin{theorem}
There exists a domain $\Omega \subset \C$ and a complete proper minimal immersion $\psi:\Omega \rightarrow \R^3$ which has uncountably many ends.
\end{theorem}
\begin{proof}
The required immersion will be obtain as a limit of a sequence of complete proper minimal immersions defined on subdomains of the complex plane. Along this section, given $\alpha$ a Jordan curve in $\C$, we denote by $\mathtt{I}(\alpha)$ as the bounded connected component of $\C-\alpha.$ In the following, we construct a sequence
$$ \chi _n=\{\Gamma_n,M_n,X_n,d_n,\ep_n, r_n\} \;, \quad \text{where}$$
\begin{enumerate}[{(a)}]
\item $\{d_n\}_{n }$, $\{\ep_n\}_{n}$ and $\{r_n\}_{n}$ are sequences of positive real numbers decreasing to zero such that $\sum_{i=n}^\infty \ep_i<r_n$.
\item \label{gamman} $\Gamma_n=\{\beta_n\}\cup\{\gamma_n(k_1,\ldots,k_j) \mid  k_i \in \{0,1\} \; , 1 \leq i \leq j \; , 1 \leq j\leq n \}\subset \C$ is a family of curves with
\begin{enumerate}[{(b}. 1)]
\item \label{b1} $\beta_n \subset \C-\D$, $n \in \N$,  are Jordan curves  wich satisfy $\overline{\mathtt{I}(\beta_{n-1})}\subset \mathtt{I}(\beta_{n})$. 
\item $\gamma_n(k_1,\ldots,k_j)$ are cycles in $\mathtt{I}(\beta_n)$ satisfying the following properties
\begin{enumerate}[{(b. 2}. 1)]
\item \label{b21} For each $1\leq j<n$ the cycle $\gamma_n(k_1,\ldots,k_j)$ is in the homology class of $\gamma_{n-1}(k_1,\ldots,k_j)$ and $\overline{\mathtt{I}(\gamma_n(k_1,\ldots,k_j))}\subset \mathtt{I}(\gamma_{n-1}(k_1,\ldots,k_j))$. 
\item  $\gamma_n(k_1,\ldots,k_n)$ is a circle centered at $c(k_1,\ldots,k_n)$ and radius $r_n>0$, where $c(k_1,\ldots,k_j)=c(k_1,\ldots,k_{j-1})+ (-1)^{k_j}\varepsilon(j) \rho(k_1,\ldots,k_j)$, with $\varepsilon(j)={\rm i}$ if $j$ is odd, $\varepsilon(j)=1$ if $j$ is even and $\rho(k_1,\ldots,k_j)$ are positive real number sufficiently small so that 
\begin{equation}\label{emparedado}
\overline{\D(c(k_1,\ldots,k_n),r_n)}\subset \D(c(k_1,\ldots,k_{n-1}),r_{n-1})-\overline{\mathtt{I}(\gamma_n(k_1,\ldots,k_{n-1}))}\;.
\end{equation}
\end{enumerate}
\end{enumerate}
Observe that $\Gamma_n$ is a multicycle in $\overline{\C}=\C \cup \{\infty\}$ with $\Sigma_{k=0}^n 2^k$ cycles. The disk $\intc(\beta_n)$ coincides with $\overline{\C}-\overline{\mathtt{I}(\beta_n)}$ and $\intc(\gamma_{n}(k_1,\ldots,k_j))=\mathtt{I}(\gamma_{n}(k_1,\ldots,k_j))$. Note also that $\overline{M(\Gamma_{n-1})}\subset M(\Gamma_n)$.
\item $M_n$ is a domain in $\C$ topologically equivalent to $M(\Gamma_n)$ and such that $\overline{M(\Gamma_n)}\subset M_n \subset M(\Gamma_n^{-\ep_n})$.

\item \label{immersionn} $X_n:M_n \longrightarrow \R^3$ is a conformal complete proper minimal immersion. 
\end{enumerate}
\begin{figure}[htbp]
	\begin{center}
		\includegraphics[width=0.80\textwidth]{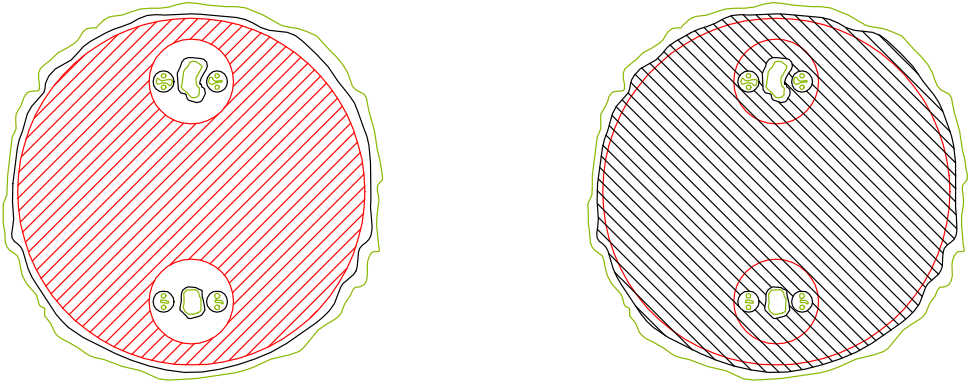}
	\end{center}
	\caption{The domains $M(\Gamma_1)$ and $M(\Gamma_2)$.}
	\label{fig:dominio}
\end{figure}
Let see that the sequence $\{\chi _n\}_{n \in \N}$ can be constructed to satisfy the following conditions
\begin{enumerate}[({Z}.$1_n$)]
\item \label{cerca} $\|X_n-X_{n-1}\|<\ep_n\;$ in $M(\Gamma_{n})\; , n\geq 2$;
\item \label{fuerabola}$X_n(M_n-M(\Gamma_{n}))\subset \R^3 -B(0,n-2 d_n-\ep_n)$, where $B(0,r)=\{x\in \R^3 \mid \|x\|<r\}$ for $r>0$;
\item \label{borde} $X_n \left( \Gamma_n \right)\subset B(0,n+\ep_n)-\overline{B(0,n-d_{n}-\ep_n)}\; .$
\item \label{landas} $\lambda_{X_n}\geq \alpha_n \:\lambda_{X_{n-1}}$ in $M(\Gamma_n)$ where $\mathcal{S}_{X_n}=\lambda_{X_n} \cdot <\cdot,\cdot> \,$, and $\{\a_n\}_{n \in \N}$ is the sequence of Remark \ref{alfas}.
\end{enumerate}
First we present the first term of the sequence. Let $X:\D\to B(0,1)$ be the immersion given by the inclusion,  $\{d_n\}_{n }$ and $\{\ep_n\}_{n \in \N}$ two sequences of positive real numbers decreasing to zero such that $\sum_{i=1}^\infty \ep_i<\frac{1}{32}$. Now we consider $\beta_1$ the circle of radius $1-\frac{1}{32}$ and center $0$, and  $\gamma_1(0)$, $\gamma_1(1)$ the circles of radius $r_1=\frac{1}{16}$ and centers $c(0)=\frac{7}{8} {\rm i}$ and $c(1)=-\frac{7}{8}{\rm i}$, respectively. Then, we can apply Corollary \ref{gordo} to the immersion $X$, the convex domain $D=B(0,1)$, the multicycle $\Gamma_1=\{\beta_1\}\cup \{\gamma_1(0),\gamma_1(1)\}$, $d=d_1=\frac{1}{4}$ and $\ep=\ep_1$ to obtain a domain $M_1$ with $\overline{M(\Gamma_1)}\subset M_1 \subset \overline{M_1} \subset M(\Gamma_1^{-\ep_1})$ and a conformal complete proper minimal immersion $X_1:M_1 \longrightarrow \R^3$ such that $\|X-X_1\|<\ep_1$ in $M(\Gamma_1)$ and $X_1(M_1-M(\Gamma_1))\subset\R^3-B(0,\frac{1}{2}-\epsilon_1)$. From here it is easy to check that (Z.\ref{borde}$_1$) is satisfied.

Assume we have constructed $\{\chi_1,\ldots,\chi_n\}$ satisfying the corresponding definitions and properties. We will define now $\chi_{n+1}$. From Property (\ref{immersionn}) and (Z.\ref{borde}$_n$) we can assert that there exist cycles $\beta_{n+1}$ and $\gamma_{n+1}(k_1,\ldots,k_j)$  fulfilling the conditions (b.1) and (b.2.1), respectively, and such that
\begin{equation*}
X_n \left(\beta_{n+1}\cup \left(\cup_{j=1}^{n}\gamma_{n+1}(k_1,\ldots,k_j)\right)\right)\subset B(0,n+1)-\overline{B(0,n+1-d_{n+1})}\; .
\end{equation*}
Furthermore, we can find $\rho(k_1,\ldots,k_{n+1})$ and $r_{n+1}>0$ appropriate so that \eqref{emparedado} is satisfied for $n+1$ and the curves $\gamma_{n+1}(k_1,\ldots,k_{n+1})$ described in (b.2.2) also fulfill the above equation, it is to say 
\begin{equation}\label{hipotesisgordo}
X_n(\Gamma_{n+1})\subset B(0,n+1)-\overline{B(0,n+1-d_{n+1})}\; .
\end{equation}
Recall that we have a sequence $\{\ep_i\}_{i}$ such that $\sum_{i=j}^\infty \ep_i<r_j$ for $j\leq n$. If $\sum_{i=n+1}^\infty \ep_i<r_{n+1}$ we do not modify the sequence $\{\ep_i\}_{i}$. If it is not the case we consider a new sequence $\{\ep_i'\}_{i}$ defined as $\ep_i'=\ep_i$ for $i \leq n$ and $\ep_i'=\ep_i r_{n+1}$ for $i \geq n+1$. It is clear that $\sum_{i=n+1}^\infty \ep'_i<r_{n+1}$. Moreover, for $j \leq n$ we have $\sum_{i=j}^\infty \ep'_i \leq \sum_{i=j}^\infty \ep_i<r_j$. For the sake of simplicity, we continue denoting  the new sequence as $\{\ep_i\}_{i}$. 

Taking into account \eqref{hipotesisgordo} we can apply Corollary \ref{gordo} to the immersion $X_n$, the multicycle $\Gamma_{n+1}$, the convex $D=B(0,n+1)$, $d=d_{n+1}$ and $\ep=\ep_{n+1}$ to obtain a domain $M_{n+1}$ with $\overline{M(\Gamma_{n+1})}\subset M_{n+1} \subset M(\Gamma_{n+1}^{-\ep_{n+1}})$ and a conformal complete proper minimal immersion $X_{n+1}:M_{n+1} \longrightarrow \R^3$ satisfying (Z.\ref{cerca}$_{n+1}$) and (Z.\ref{fuerabola}$_{n+1}$). Moreover, from   (Z.\ref{cerca}$_{n+1}$) we obtain that if $\ep_{n+1}$ is sufficiently small the property (Z.\ref{landas}$_{n+1}$) is also satisfied. Finally, (Z.\ref{borde}$_{n+1}$) follows from \eqref{hipotesisgordo} and (Z.\ref{cerca}$_{n+1}$). Consequently, we have the sequence $\{\chi_n\}$. 

Hereafter, we define the required immersion $\psi$. We denote by $\Omega=\bigcup_{n=1}^{\infty} M(\Gamma_n)$. Clearly $\Omega$ is a domain, since it is the union of domains with non empty intersection. Furthermore, from  (Z.\ref{cerca}$_{n}$) we deduce that the sequence $\{X_n\}_{n\in \N}$ converges uniformly on compact sets of $\Omega$ and so we can define $\psi:\Omega  \to \R^3$ as $\psi(z)=\lim_{i \to  \infty,i \geq n}X_i(z)$ for $z \in M(\Gamma_n)$. By making use of Harnack's theorem we know that $\psi$ is a harmonic map. Let us see that $\psi$ is immersion. Take $z \in \Omega$. Thus, there exists $n \in \N$ such that $z \in M(\Gamma_n)$. Then, according to properties (Z.\ref{landas}$_{m}$) for  $m \geq n$, we have 
$$ \lambda_{X_m}\geq \alpha_m \:\lambda_{X_{m-1}} \geq \cdots \geq \prod_{i=n}^m \alpha_i \, \lambda_{X_{n-1}} \geq \frac{1}{2} \, \lambda_{X_{n-1}}>0 \; .$$
By taking limits in the above inequality as $m \to \infty$ we obtain $\lambda_\psi(z)>0$. On the other hand, it is easy to obtain the properness (and therefore the completeness) of $\psi$ from properties (Z.\ref{fuerabola}$_{n}$). 

Finally, let us demonstrate that $\psi:\Omega \rightarrow \R^3$ possesses uncountably many ends. Let $Q$ denote a sequence $Q=\{k_i\}_{i \in \N}$, where $k_i \in \{0,1\}$. Next, we consider any proper arc $\sigma_Q:[0,\infty[ \; \to \Omega$ satisfying 
\begin{equation}\label{arcos}
\sigma_Q([j,\infty[) \subset \D(c(k_1,\dotsc, k_j),r_j) \; , \forall \; j \in \N\; . 
\end{equation}
We note first that if $Q=\{k_i\}_{i \in \N}$ and $Q'=\{k_i'\}_{i \in \N}$ are two sequences as above such that $Q \neq Q'$ then there exists $j_0=\min\{j \in \N \mid k_j \neq k_j'\}$. Thus, \eqref{arcos} implies that $\D(c(k_1,\dotsc, k_{j_0}),r_{j_0})$ and $\D(c(k_1',\dotsc, k_{j_0}'),r_{j_0})$ are two disks containing $\sigma_Q([j_0,\infty[)$ and $\sigma_{Q'}([j_0,\infty[)$, respectively. Since 
$$\overline{\D(c(k_1,\dotsc, k_{j_0}),r_{j_0})} \cap \overline{\D(c(k_1',\dotsc, k_{j_0}'),r_{j_0})}=\emptyset \;,$$ we can consider $\partial(\D(c(k_1,\dotsc, k_{j_0}),r_{j_0}))$ as a compact set separating $\sigma_Q([j_0,\infty[)$ and $\sigma_{Q'}([j_0,\infty[)$. Therefore, $\sigma_Q$ and $\sigma_{Q'}$ are two distinct topological ends. As there exists an uncountable number of sequences $Q$, we deduce that there are uncountably many ends.
\end{proof}

Finally, we would like to mention the following:
\begin{remark}
With the same ideas presented in the proof of the above theorem it is also possible to construct properly immersed minimal surfaces with uncountably many ends in such a way that all the ends are limit ends. 
\end{remark}




\end{document}